\documentclass[12pt]{article}
\usepackage{amsmath,amssymb,amsfonts,amscd,amsxtra}
\usepackage{latexsym}

\input{epsf.sty}
\usepackage{epsfig}
\usepackage{graphicx,color}
\usepackage{graphics}
\usepackage{subfigure}
\usepackage{appendix}
\usepackage[top=1in, bottom=1in, left=1in, right=1in]{geometry}

\renewcommand{\theequation}{\thesection.\arabic{equation}}

\newtheorem{thm}{Theorem}[section]

\newtheorem{lem}[thm]{Lemma}
\newtheorem{prop}[thm]{Proposition}

\newtheorem{rmk}[thm]{Remark}

\renewcommand{\Re}{{\mbox{Re}}}

\newcommand{\abs}[1]{\left\vert#1\right\vert}

\newcommand{\qed}{\hfill \ensuremath{\square}}

\renewcommand\appendix{\par
  \setcounter{section}{0}
  \setcounter{subsection}{0}
  \setcounter{figure}{0}
  \setcounter{table}{0}
  \renewcommand\thesection{Appendix \Alph{section}}
  \renewcommand\theequation{\Alph{section}.\arabic{equation}}
  \renewcommand\thefigure{\Alph{section}.\arabic{figure}}
  \renewcommand\thetable{\Alph{section}.\arabic{table}}
  \renewcommand\thethm{\Alph{section}.\arabic{thm}}
}

\numberwithin{equation}{section}

\date{}

\title{Scattering by a periodic array of subwavelength slits II: surface bound state, total transmission and field enhancement in homogenization regimes}
\author{
Junshan Lin \thanks{\footnotesize Department of Mathematics and Statistics, Auburn University, Auburn, AL 36849 (jzl0097@
auburn.edu). Junshan Lin was partially supported by the NSF grant DMS-1417676.}
 \; and Hai Zhang\thanks{\footnotesize 
Department of Mathematics, 
 HKUST,  Clear Water Bay, Kowloon, Hong Kong (haizhang@ust.hk). Hai Zhang was supported by HK RGC grant ECS 26301016 and the UGC grant SBI17SC12 from HKUST.}}

\begin{document}

\maketitle

\begin{abstract}
This is the second part in a series of two papers that concern with the quantitative analysis of the electromagnetic field enhancement and anomalous diffraction by a periodic array of subwavelength slits. In this part, we explore the scattering problem in the homogenization regimes, where the size of the period is much smaller than the incident wavelength. In particular, two homogenization regimes are investigated, where the size of the pattered slits has the same order as the size of the period in the first configuration, and the size of the slit is much smaller than the size of the period in the second configuration.
By presenting rigorous asymptotic analysis, we demonstrate that surface plasmonic effect mimicking that of plasmonic metals occurs in the first homogenization regime. The corresponding dispersion curve lies below the light line and the associated eigenmodes are surface bound sates.
In addition, for the incident plane wave, we discover and justify a novel phenomenon of total transmission which occurs either at certain frequencies for all incident angles, or at a special incident angle but for all frequencies. For the second homogenization regime, the non-resonant field enhancement is investigated, and it is shown that the fast transition of the magnetic field in the slit induces strong electric field enhancement. Moreover, the enhancement becomes stronger when the coupling of the slits is weaker. 
 \end{abstract}


\setcounter{equation}{0}
\setlength{\arraycolsep}{0.25em}
\section{Introduction}
\begin{figure}[!htbp]
\begin{center}
\includegraphics[height=4.2cm,width=15cm]{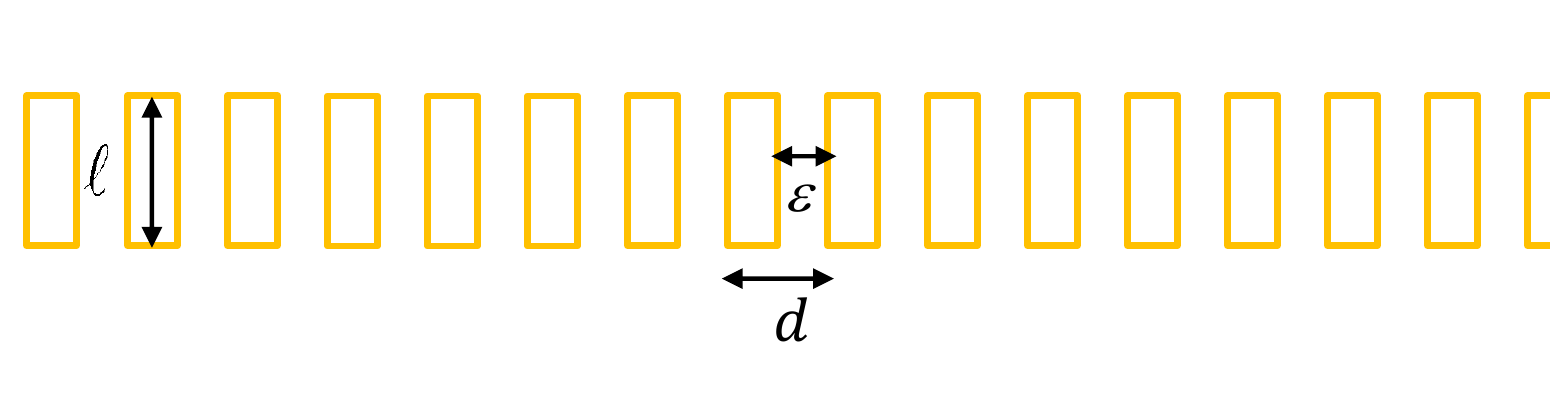}\\
\vspace*{-15pt}
\includegraphics[height=4.2cm,width=15cm]{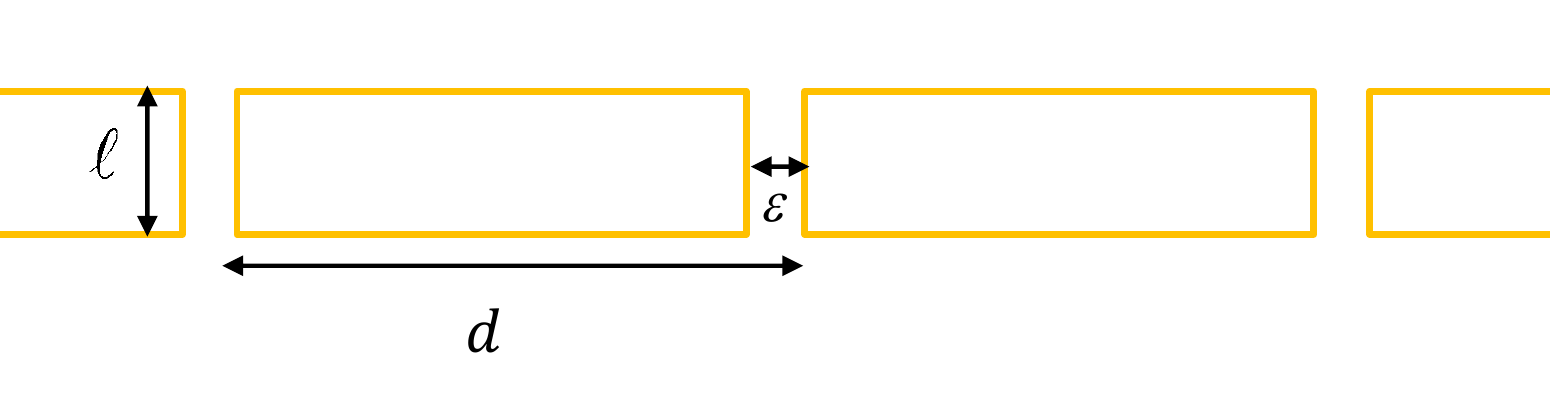}
\caption{Geometry of slits in two homogenization regimes.
Top:  $\ell =1$, $\varepsilon \sim d \ll 1$; Bottom: $\ell =1$, $\varepsilon \ll 1$,  $d\sim 1$, 
or $d\gg1$ but $d\ll\lambda$.}
\label{fig-geo_H1H2}
\end{center}
\end{figure}

This is the second part in a series of two papers that are concerned with the electromagnetic scattering and
field enhancement for a perfect conducting (PEC) slab patterned with a periodic array of subwavelength slits.
In the first part \cite{lin_zhang17}, we investigated the field enhancement in the diffraction regime, where the size of the period
is the same order as the incident wavelength. In this paper, we explore the scattering problem in the homogenization regime,
where the size of the period is much smaller than the incident wavelength. Similar settings have been investigated for periodically arranged subwavelength resonators such as plasmonic particles and bubbles in \cite{hai11, hai22}, where the mechanism of metasurface is explained. 
We shall consider two homogenization regimes. In the first regime,  the size of the pattered slits has the same order as the size of the period 
(see Figure \ref{fig-geo_H1H2}, top),
while in the second regime, the size of the slit is much smaller than the size of period (see Figure \ref{fig-geo_H1H2}, bottom).
The studies are motiavted by recent growing interest in extraordinary optical transmission and strongly enhanced electromagnetic fields
in subwavelength apertures or holes, which could lead to potentially significant applications in biological and chemical sensing, near-field spectroscopy, etc 
\cite{chen14, ebbesen98, garcia05, garcia10, krie04, pendry04}. 
The readers are also referred  to \cite{lin_zhang16} for scattering and field enhancement for a single narrow slit and
\cite{eric10, eric10-2} for a closely related problem of scattering by subwavelength cavities.

\begin{figure}[!htbp]
\begin{center}
\includegraphics[height=5.5cm,width=15cm]{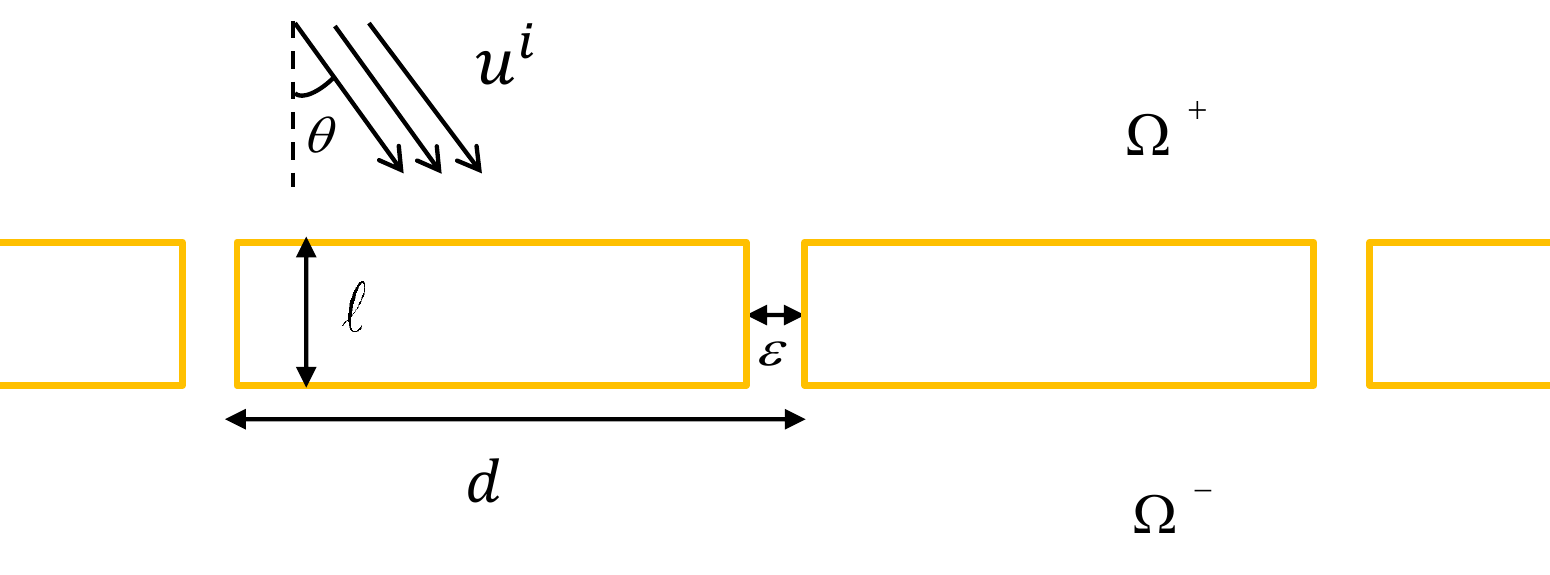}\\
\caption{Setup of the scattering problem. 
The slits $S_\varepsilon$ are arranged periodically with the size of the period $d$, and each slit has a rectangular shape of length $\ell$ and width $\varepsilon$ respectively. The domains above and below the perfect conductor slab are denoted as $\Omega^{+}$ and $\Omega^{-}$ respectively, and the domain exterior to the perfect conductor is denoted as $\Omega_{\varepsilon}$, which consists of $S_\varepsilon$, $\Omega^{+}$, and $\Omega^{-}$. }\label{fig-prob_geo}
\end{center}
\end{figure}

We now present the setup of the scattering problem. Figure \ref{fig-prob_geo} depicts the geometry of the cross section for the metallic structure under consideration. 
The slab occupies the domain 
$\{(x_1,x_2)\;|\; 0<x_2<\ell\}$ on the $x_1x_2$ plane, where $\ell$ is the thickness of the metal slab. The slits, which are invariant along the $x_3$ direction, 
occupy the region $\displaystyle{S_\varepsilon=\bigcup_{n=0}^{\infty} S_\varepsilon^{(0)} + nd}$, where $d$ is the size of the period, and
$S_\varepsilon^{(0)}:=\{(x_1,x_2)\;|\; 0<x_1<\varepsilon, 0<x_2<\ell \}$ is of rectangular shape.
We denote the semi-infinite domain above and below the slab by $\Omega^{+}$ and $\Omega^{-}$,
and $\Omega_{\varepsilon}$ the domain exterior to the perfect conductor, i.e.,
$\Omega_{\varepsilon}=\Omega^{+}\cup\Omega^{-}\cup S_\varepsilon$. We also denote by $\nu$  the unit outward normal pointing to the exterior domain $\Omega^{+}$ or $\Omega^{-}$.

The width of slit, $\varepsilon$, is assumed to be much smaller than the thickness of the slab $\ell$.
For clarity of exposition, we shall set $\ell=1$ in all technical derivations. The general case
for $\ell\neq1$ follows by a normalization process and a scaling argument. Furthermore, we assume that the size of the period $d$
is much smaller than the wavelength $\lambda$ such that the problem under consideration is in the homogenization regime.
The following two homogenization regimes are investigated here:
\begin{itemize}
\item [(H1)] The scaling of geometrical parameters are given by $\ell =1$,
$\varepsilon \sim d \ll 1$, and the incident wavelength $\lambda \sim O(1)$ or $\lambda \gg 1$. That is, $\varepsilon \sim d \ll \lambda$.
A schematic plot of the geometry is shown in Figure \ref{fig-geo_H1H2} (top).
\item [(H2)] The scaling of geometrical parameters are given by $\ell =1$, $\varepsilon \ll 1$,  $d\sim 1$ or  $1 \ll d \ll \lambda$, 
and $\lambda \gg 1$. That is $\varepsilon \ll d \ll \lambda$. A schematic plot of the geometry is shown in Figure \ref{fig-geo_H1H2} (bottom).
\end{itemize}

Assume that a polarized time-harmonic electromagnetic wave impinges upon the perfect conductor from the above.
We consider the transverse magnetic (TM) case where the incident magnetic field is perpendicular to the  $x_1x_2$ plane,
and its $x_3$ component is given by the scalar function $u^{i}= e^{i (\kappa x_1 - \zeta (x_2-1))}$. Here
$\kappa= k \sin \theta$, $\zeta=k\cos\theta$,  $k$ is the wavenumber and $\theta$ is the incident angle. 
Throughout the paper, we assume that $|\theta| < \theta_0<\dfrac{\pi}{2}$ for some $\theta_0$ to exclude the case of grazing incidence angle.
The total field $u_{\varepsilon}$, which consists of the incident wave $u^{i}$ and the scattered field $u_{\varepsilon}^s$, satisfies the Helmholtz equation
\begin{equation}\label{eq-Helmholtz}
\Delta u_{\varepsilon} + k^2 u_{\varepsilon} = 0   \quad\quad  \mbox{in} \; \Omega_{\varepsilon},
\end{equation}
and the boundary conditition
\begin{equation}\label{eq-bnd_cnd}
\dfrac{\partial u_{\varepsilon}}{\partial \nu} = 0  \quad \mbox{on} \; \partial \Omega_{\varepsilon}.
\end{equation}
 We look for quasi-periodic solutions such that $u_{\varepsilon}(x_1,x_2)=e^{i\kappa x_1}\tilde u_{\varepsilon}(x_1,x_2)$,
where $\tilde u_{\varepsilon}$ is a periodic function with $\tilde u_{\varepsilon}(x_1+d,x_2)=\tilde u_{\varepsilon}(x_1,x_2)$,
or equivalently,
\begin{equation}\label{eq-quasi_periodic}
u_{\varepsilon}(x_1+d,x_2)=e^{i\kappa d} u_{\varepsilon}(x_1,x_2).
\end{equation}
Define
$$
\kappa_n=\kappa+\dfrac{2\pi n}{d} \quad \mbox{and} \quad
\zeta_n(k)= \sqrt{k^2-\kappa_n^2},
$$
where the function $f(z)=\sqrt{z}$ is understood as an analytic function defined in the domain $\mathbf{C}\backslash\{-it: t\geq 0\}$ by
$$
\sqrt{z} = |z|^{\frac{1}{2}} e^{\frac{1}{2} i\arg z}
$$
throughout the paper.  Then it can be shown that the outgoing scattered field adopts the following Rayleigh-Bloch expansion in
$\Omega^{+}$ and $\Omega^{-}$ respectively (cf. \cite{bao95, bonnet_starling-94, shipman-10})
\begin{equation}\label{eq-rad_cond}
 u_{\varepsilon}^s(x_1,x_2) =  \sum_{n=-\infty}^{\infty} u_n^{s,+} e^{i \kappa_n x_1 + i\zeta_n  x_2 } \quad \mbox{and} \quad
 u_{\varepsilon}^s(x_1,x_2) =  \sum_{n=-\infty}^{\infty} u_n^{s,-} e^{i \kappa_n x_1 - i\zeta_n  x_2 },
\end{equation}
where $u_n^{s,\pm}$ are constants.  
The expansion \eqref{eq-rad_cond} is usually referred to as the outgoing radiation condition and is imposed for the scattered field in the semi-infinite domains. 
In sum, the mathematical model for the scattering problem is defined in the domain $\Omega_{\varepsilon}$ and
given by the equations \eqref{eq-Helmholtz} - \eqref{eq-rad_cond}.
Due to the quasi-periodicity of the solution, we will restrict $\kappa$ to the first Brillouin zone  $(-\pi/d,\pi/d]$.
Such $\kappa$ is called the reduced wave vector component \cite{bonnet_starling-94, shipman-10}.

In this paper, based upon a combination of layer potential techniques and asymptotic analysis, 
we develop a quantitative analysis of field enhancement and anomalous transmission behavior for the scattering problem
in the above mentioned two homogenization regimes. In more details,
\begin{itemize}
\item [(i)]  In the homogenization regime (H1), the asymptotic expansions of the dispersion relation and the associated eigenmodes are derived.
It is demonstrated that surface plasmonic effect mimicking that of plasmonic metals occurs in such a configuration. 
More precisely, the dispersion curve, which lies below the light line with $k(\kappa)<|\kappa|$, resembles that of surface plasmon polaritons of the nobel metal slab; and the eigenmodes, which are surface bound states along the boundaries the perfect conducting slab, resemble the plasmonic waves of noble metals.
Therefore, the specific configuration with $\varepsilon \sim d \ll \lambda$ in this homogenization regime extends the frequency band for the surface plasmon, which is originally supported on noble metals in optical and near-infrared regime, to the lower frequency regime where metals can be viewed as perfect conductors. This is so-called spoof surface plasmon in physics literatures, and has the potential for openning
new opportunities to control radiation at surfaces over a wide spectral range \cite{garcia05, pendry04}.

\item [(ii)] We derive the asymptotic expansions of the scattered wave field when an incident plane wave impinges on the periodic structure as specified in the homogenization regime (H1). In such a scenario, the reduce wave vector component satisfies $|\kappa|=|k\sin\theta|<k$ and the solution to the scattering problem  \eqref{eq-Helmholtz} - \eqref{eq-rad_cond}  is unique. Interestingly, it is shown that total transmission through small slits can be achieved either at certain frequencies for all incident angles or for all frequencies at a specific incident angle. We clarify that such perfect transmission is not due to plasmonic resonant effect or scattering resonance. Instead, for the former, it might be related to Fabry-Perot resonance associated with the homogenized homogeneous slab, where all reflected waves from the slab boundaries interfere destructively \cite{vaughan89}. The mechanism for the latter is not clear to us. 

\item[(iii)] In the homogenization regime (H2), there exists no complex resonance or real eigenvalue,
and the scattering problem  \eqref{eq-Helmholtz} - \eqref{eq-rad_cond} attains a unique solution.
We derive the asymptotic expansion of the electromagnetic fields in the near and far field, and show that although no enhancement is gained for the magnetic field, strong electric field is induced in the slits and on the slit apertures. Such field enhancement is not induced by resonances, but due to the fast transition of the magnetic field in the slits. In addition, we also discuss the enhancement behavior with varying sizes of the period $d$. We show that as the period $d$ decreases and the coupling between the slits is stronger, the field enhancement becomes weaker. 
\end{itemize}

The rest of the paper is organized as follows. We begin by introducing layer potentials for the scattering problem 
and presenting the asymptotic expansion for the solution to the scattering problem in Section \ref{sec-bnd_int_eqn} for both homogenization regimes. 
The quantitative analysis of anomalous transmission and field enhancement behaviors
is presented in Section \ref{sec-H1} and \ref{sec-H2} for the homogenization regime (H1) and (H2) respectively.
The paper is concluded with some remarks about the ongoing and future works along this direction in Section \ref{sec-conclusion}.

\setcounter{equation}{0}
\section{Boundary integral equations and the solution to the scattering problem}\label{sec-bnd_int_eqn}
\subsection{Layer potentials and boundary integral formulations}
In this section, we collect some preliminaries on the layer potentials and boundary integral formulations for the scattering problem.
The readers are referred to the first part of this series \cite{lin_zhang17} for the proof.
For a given $\kappa\in(-\pi/d,\pi/d]$, let 
\begin{equation}\label{eq-gd}
g^d(x,y)=g^d(x,y; \kappa) = -\dfrac{i}{2d}  \sum_{n=-\infty}^{\infty} \dfrac{1}{\zeta_n(k)   } e^{i \kappa_n(x_1-y_1)+i\zeta_n(k)   |x_2-y_2| },
\end{equation}
where
\begin{equation*}
\kappa_n=\kappa+\dfrac{2\pi n}{d} \quad \mbox{and} \quad
\zeta_n(k)   = \left\{
\begin{array}{lll}
\vspace*{5pt}
\sqrt{k^2-\kappa_n^2},  & \abs{\kappa_n} < k, \\
i\sqrt{\kappa_n^2-k^2},  & \abs{\kappa_n} > k. \\
\end{array}
\right.
\end{equation*}
It is clear that $g^d(x,y, \kappa)$ is the periodic Green's function which solves the following equation
$$ \Delta g^d(x,y;\kappa) + k^2 g^d(x,y;\kappa)=e^{i\kappa(x_1-y_1)} \sum_{n=-\infty}^{\infty} \delta(x_1-y_1-nd)\delta(x_2-y_2) \quad x, y\in\mathbf{R}^2.  $$

The exterior Green's function $g^e(x,y)=g^e(x,y, \kappa)$ in domain $\Omega^+\cup\Omega^-$ with the Neumann boundary condition $\dfrac{\partial g^e(x,y, \kappa)}{\partial \nu_y}=0$ on $\{y_2=1\}$ and $\{y_2=0\}$ is then given by $g^e(x,y) =g^d(x,y, \kappa)+g^d(x',y, \kappa)$, where 
\begin{equation*}
x' = \left\{
\begin{array}{ll}
(x_1, 2-x_2) & \mbox{if} \; x, y \in \Omega^+,  \\
(x_1,-x_2) &  \mbox{if} \; x, y \in \Omega^-. 
\end{array}
\right.
\end{equation*}

The Green's function $g_\varepsilon^i(x,y)$ that solves
$$ \Delta g_\varepsilon^i(x,y, k) + k^2 g_\varepsilon^i(x,y, k) = \delta(x-y), \quad x, y\in  S_{\varepsilon}^{(0)} $$
with the Neumann boundary condition may be expressed as
$$ g_\varepsilon^i(x,y)= \sum_{m,n=0}^\infty c_{mn}\phi_{mn}(x)\phi_{mn}(y), $$
where $c_{mn}=\dfrac{1}{k^2-(m\pi/\varepsilon)^2 - (n\pi)^2}$,  $\phi_{mn}=\sqrt{\dfrac{a_{mn}}{\varepsilon}}\cos\left(\dfrac{m\pi x_1}{\varepsilon}\right) \cos(n\pi x_2)$
with the coefficient
\begin{equation*}
a_{mn} = \left\{
\begin{array}{llll}
1  & m=n=0, \\
2  & m=0, n\ge 1 \quad \mbox{or} \quad n=0, m\ge 1, \\
4  & m\ge 1, n \ge 1.
\end{array}
\right.
\end{equation*}

\begin{figure}[!htbp]
\begin{center}
\includegraphics[height=5.6cm,width=14cm]{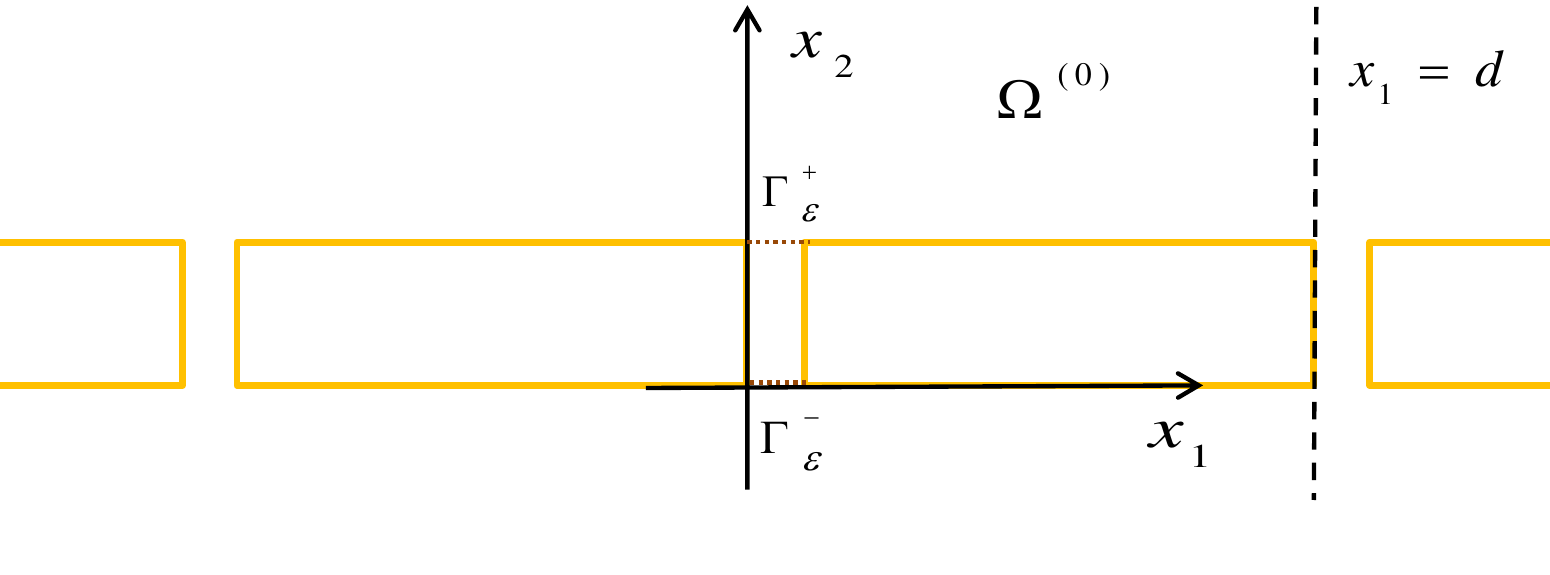}
\vspace{-10pt}
\caption{Problem geometry in one reference period $\Omega^{(0)}$.}\label{fig-prob_one_period}
\vspace{-10pt}
\end{center}
\end{figure}

To define the layer potentials, we consider the reference period $\Omega^{(0)}:=\{ x \in \mathbf{R}^2 \; | \;   0<x_1<d \}$ 
as shown in Figure \ref{fig-prob_one_period}.
Denote the  the upper and lower aperture of the slit $S_\varepsilon^{(0)}$ in the reference period $\Omega^{(0)}$ 
by $\Gamma^{+}_\varepsilon$ and $\Gamma^{-}_\varepsilon$ respectively (see Figure \ref{fig-prob_one_period}).
\begin{lem}[\cite{lin_zhang17}]\label{eq-u_eps_formula}
Let $u_\varepsilon(x)$ be the solution of the scattering problem \eqref{eq-Helmholtz} - \eqref{eq-rad_cond}, then
\begin{eqnarray*}
u_\varepsilon(x) &=& \int_{\Gamma^+_\varepsilon} g^e(x,y) \dfrac{\partial u_\varepsilon(y)}{\partial y_2} ds_y + u^i+ u^r  \quad \mbox{for} \; x\in\Omega^{(0)} \cap \Omega^+. \\
u_\varepsilon(x) &=& -\int_{\Gamma^-_\varepsilon} g^e(x,y) \dfrac{\partial u_\varepsilon(y)}{\partial y_2} ds_y  \quad \mbox{for} \; x\in\Omega^{(0)} \cap \Omega^-. \\
u_\varepsilon(x) &=& \int_{\Gamma^-_\varepsilon} g_\varepsilon^i(x,y) \dfrac{\partial u_\varepsilon(y)}{\partial y_2} ds_y-\int_{\Gamma^+_\varepsilon} g_\varepsilon^i(x,y) \dfrac{\partial u_\varepsilon(y)}{\partial y_2} ds_y \quad \mbox{for} \; x\in S_\varepsilon^{(0)}.
\end{eqnarray*}
Here $u^r=e^{i (\kappa x_1 + \zeta (x_2-1))}$ is the reflected field of the ground plane $\{x_2=1\}$ without slits.
\end{lem}

Based upon Lemma  \ref{eq-u_eps_formula} and the continuity of the single layer potential, we obtain the following boundary integral equations
defined over the slit apertures $\Gamma^\pm_\varepsilon$.
\begin{lem} \label{lem-ppp}
The following hold for the solution of the scattering problem \eqref{eq-Helmholtz} - \eqref{eq-rad_cond}:
\begin{equation}\label{eq-T1_e_eps}
u_\varepsilon(x) = \int_{\Gamma^+_\varepsilon}  g^e(x,y)\dfrac{\partial u_\varepsilon(y)}{\partial y_2} ds_y + u^i+u^r \quad \mbox{for} \; x\in\Gamma^+_\varepsilon.
\end{equation}
\begin{equation}\label{eq-T2_e_eps}
u_\varepsilon(x) = -\int_{\Gamma^-_\varepsilon}  g^e(x,y) \dfrac{\partial u_\varepsilon(y)}{\partial y_2} ds_y \quad \mbox{for} \; x\in\Gamma^-_\varepsilon.
\end{equation}
\begin{equation}\label{eq-T_i_eps}
u_\varepsilon(x) =  \int_{\Gamma^-_\varepsilon} g_\varepsilon^i(x,y) \dfrac{\partial u_\varepsilon(y)}{\partial y_2} ds_y-\int_{\Gamma^+_\varepsilon} g_\varepsilon^i(x,y) \dfrac{\partial u_\varepsilon(y)}{\partial y_2} ds_y
 \quad \mbox{for} \; x\in \Gamma^+_\varepsilon\cup\Gamma^-_\varepsilon.
\end{equation}
\end{lem}
An application of the above Lemma leads to the following system of integral equations:
\begin{equation} \label{eq-scattering2}
\left\{
\begin{array}{lll}
&\displaystyle{\int_{\Gamma^+_\varepsilon}  g^e(x,y)\dfrac{\partial u_\varepsilon(y)}{\partial y_2} ds_y 
+\int_{\Gamma^+_\varepsilon} g_\varepsilon^i(x,y) \dfrac{\partial u_\varepsilon(y)}{\partial y_2} ds_y
-\int_{\Gamma^-_\varepsilon} g_\varepsilon^i(x,y) \dfrac{\partial u_\varepsilon(y)}{\partial y_2} ds_y
+ u^i+u^r =0}, \quad \mbox{on} \,\, \Gamma^+_\varepsilon, \\ \\
& -\displaystyle{\int_{\Gamma^-_\varepsilon}  g^e(x,y)\dfrac{\partial u_\varepsilon(y)}{\partial y_2} ds_y 
+\int_{\Gamma^+_\varepsilon} g_\varepsilon^i(x,y) \dfrac{\partial u_\varepsilon(y)}{\partial y_2} ds_y
-\int_{\Gamma^-_\varepsilon} g_\varepsilon^i(x,y) \dfrac{\partial u_\varepsilon(y)}{\partial y_2} ds_y=0,} \quad \mbox{on} \,\, \Gamma^-_\varepsilon. 
\end{array}
\right.
\end{equation}

\begin{prop}
The scattering problem \eqref{eq-Helmholtz} - \eqref{eq-rad_cond} is equivalent to the system of boundary integral equations (\ref{eq-scattering2}). 
\end{prop}

It is clear that 
$$
\left.\dfrac{\partial u_\varepsilon}{\partial \nu}\right|_{\Gamma^+_\varepsilon} =\dfrac{\partial u_\varepsilon}{\partial y_2}(y_1, 1), 
\quad \left.\dfrac{\partial u_\varepsilon}{\partial \nu}\right|_{\Gamma^-_\varepsilon} =-\dfrac{\partial u_\varepsilon}{\partial y_2}(y_1, 0),\quad
(u^i+u^r)|_{\Gamma^+_\varepsilon}= 2e^{i\kappa x_1}.
$$
The above functions are defined over narrow slit apertures with size $\varepsilon \ll 1$. We rescale the functions by introducing
$X=x_1/\varepsilon $ and $Y=y_1/\varepsilon$, and define the following quantities: 
\begin{eqnarray*}
&& \varphi_1(Y):=  -\dfrac{\partial u_\varepsilon}{\partial y_2}( \varepsilon Y, 1); \\
&& \varphi_2(Y):=  \dfrac{\partial u_\varepsilon}{\partial y_2}(\varepsilon Y, 0);\\ 
&& f(X):= (u^i+u^r)(\varepsilon X, 1) = 2e^{i \kappa \varepsilon X};\\
&& G_\varepsilon^e(X, Y)=G_\varepsilon^e(X, Y, \kappa) :=  g^e(\varepsilon X,1; \varepsilon Y,1)=  g^e(\varepsilon X,0; \varepsilon Y,0)=-\dfrac{i}{d} \sum_{n=-\infty}^{\infty} \dfrac{1}{\zeta_n(k)   } e^{i \kappa_n\varepsilon(X-Y)}; \\
&& G_\varepsilon^i(X, Y) := g_\varepsilon^i(\varepsilon X, 1; \varepsilon Y, 1 ) = g_\varepsilon^i(\varepsilon X, 0; \varepsilon Y, 0 ) = \sum_{m,n=0}^\infty\dfrac{c_{mn}a_{mn}}{\varepsilon}\cos(m\pi X) \cos(m\pi Y); \\
&& \tilde G_\varepsilon^i(X, Y) := g_\varepsilon^i(\varepsilon X, 1; \varepsilon  Y, 0 ) = g_\varepsilon^i(\varepsilon  X, 0; \varepsilon  Y, 1 ) =\sum_{m,n=0}^\infty\dfrac{ (-1)^n c_{mn} a_{mn}}{\varepsilon}\cos(m\pi X) \cos(m\pi Y); 
\end{eqnarray*}
We also define three boundary integral operators:
\begin{eqnarray}\label{op_T}
  && (T^e \varphi)(X) = \int_0^1 G_\varepsilon^e(X, Y)  \varphi(Y) dY  \quad  X\in (0,1); \label{op_Te} \\
 && (T^i  \varphi) (X) = \int_0^1  G_\varepsilon^i(X, Y) \varphi(Y) dY  \quad  X\in (0,1); \label{op_Ti} \\
 && (\tilde T^i  \varphi)(X)  = \int_0^1 \tilde G_\varepsilon^i(X, Y) \varphi(Y) dY  \quad  X\in (0,1).  \label{op_tilde_Ti}
\end{eqnarray}
By a change of variable $x_1=\varepsilon X$ and $y_1=\varepsilon Y$ in \eqref{eq-scattering2}, the following proposition follows.
\begin{prop}
The system of equations (\ref{eq-scattering2}) is equivalent to the following one:
\begin{equation}\label{eq-scattering3}
\left[
\begin{array}{llll}
T^e+T^i    & \tilde{T}^i \\
\tilde{T}^i  & T^e+T^i 
\end{array}
\right] \left[
\begin{array}{llll}
\varphi_1    \\
\varphi_2 
\end{array}
\right]=\left[
\begin{array}{llll}
f/\varepsilon   \\
0
\end{array}
\right].
\end{equation}
\end{prop}

\subsection{Asymptotic expansion of the boundary integral operators}
We recall several function spaces to be used throughout the paper, which are introduced in the first part \cite{lin_zhang17}.
Let $H^s(\mathbf{R})$ be the standard fractional Sobolev space for $s \in \mathbf{R}$.
For a bounded open interval $I$, define the Hilbert spaces
$$
H^s(I) := \{ u=U|_{I}  \;\big{|}\;  U \in H^s(\mathbf{R}) \}. 
$$
and
$$
\tilde H^s(I) := \{ u=U|_{I}  \;\big{|}\;  U \in H^s(\mathbf{R})\,\,  \mbox{and} \,\, supp \,U \subset \bar I   \}.
$$
Then $\tilde H^s(I)$ is the dual of $H^s(I)$.
 For simplicity of notation, we denote $V_1 = \tilde H^{-\frac{1}{2}}(0, 1)$ and 
$V_2 =  H^{\frac{1}{2}}(0, 1)$. The duality between $V_1$ and $V_2$ will be denoted by 
$\langle u, v \rangle$ for any $u \in V_1$, $v \in V_2$. 

Let us define the operator $P: V_1 \to V_2$ by 
\begin{equation}\label{opt_P}
P \varphi(X) = \langle \varphi, 1 \rangle 1,
\end{equation}
where $1$ is a function defined on the interval $(0, 1)$ and is equal to one therein. Then $1 \in V_2$ and the above definition is valid.

To obtain the solution of the scattering problem, we begin with the asymptotic expansion of the integral operators $T^e$, $T^i$ and  $\tilde T^i$.
First, the kernels $G_\varepsilon^i(X, Y)$ and $\tilde G_\varepsilon^i(X, Y)$ attain the following asymptotic expansions.
\begin{lem} \label{lem-green_slit}
Let
\begin{equation}\label{beta_i}
\beta^i(k, \varepsilon)= \dfrac{\cot k }{k \varepsilon} +  \dfrac{2\ln 2}{\pi}, \quad \tilde \beta(k,\varepsilon) = \dfrac{1}{(k\sin k) \varepsilon}, \\
\end{equation}
\begin{equation}\label{rho_i}
\rho^i(X,Y) = \dfrac{1}{\pi} \left[ \ln \left(\abs{\sin \left(\frac{\pi(X+Y)}{2}\right)}\right) + \ln \left(\abs{\sin \left(\frac{\pi(X-Y)}{2}\right)}\right) \right].
\end{equation}
If $k\varepsilon \ll 1$, then
\begin{eqnarray*}
G_\varepsilon^i(X, Y) &=& \beta^i(k ,\varepsilon) +  \rho^i(X,Y) +  r^{i}_{\varepsilon}((X, Y), \label{Gi_exp} \\ 
\tilde G_\varepsilon^i(X, Y) &= &  \tilde\beta(k,\varepsilon) + \tilde r_{\varepsilon}(X, Y). \label{tGi_exp}
\end{eqnarray*}
Here $r^{i}_{\varepsilon}(X, Y)$ and  $\tilde r(X, Y)$ are bounded functions with $r^{i}_{\varepsilon}\sim O((k\varepsilon)^2)$, and 
$\tilde r_\infty\sim O(e^{-1/\varepsilon})$ for all $X, Y\in(0,1)$.
\end{lem}
The proof of the Lemma can be found in \cite{lin_zhang16, lin_zhang17}.

\begin{lem} \label{lem-periodic_green}
Assume that $\kappa\in(-\pi/d,\pi/d]$, and $\kappa \sim O(1)$ satisfying $|\kappa/\sqrt{k^2-\kappa^2}|\le C$, where $C$ is a positive constant. 
Then the kernel $G_\varepsilon^e(X, Y)$ attains the following asymptotic expansion in both homogenization regimes (H1) and (H2):
\begin{equation}\label{Ge_exp}
G_\varepsilon^e(X, Y)= \beta^e(k,\kappa, d, \varepsilon) + \rho^e (X,Y; k,\kappa) +  r^e_\varepsilon(X,Y; k,\kappa),
\end{equation}
where $\beta^e(k,\kappa, d, \varepsilon)$ is independent of $X$ and $Y$, $\rho^e(X,Y; k,\kappa)$ is a function independent of $\varepsilon$, 
and $r^e_\varepsilon(X,Y)$ is a bounded function with $r^e_\varepsilon \sim O(r(\varepsilon))$, and $r(\varepsilon)\to 0$ as $\varepsilon\to0$.

\begin{enumerate}
\item[(1)] For the homogenization regime (H1),
\begin{equation}\label{beta_e_H1}
\beta^e(k,\kappa, d, \varepsilon)=\dfrac{1}{\pi}\ln 2 -\dfrac{i \eta}{\sqrt{k^2-\kappa^2} \; \varepsilon},
\end{equation}
\begin{equation}\label{rho_e_H1}
\rho^e (X,Y; k,\kappa) = \dfrac{1}{\pi}  \ln |\sin(\pi\eta(X-Y))| + \dfrac{\kappa \eta }{ \sqrt{k^2-\kappa^2}} (X-Y),
\end{equation}
where $\eta=\varepsilon/d$. In addition, $r(\varepsilon) = \varepsilon$ if $\kappa\neq 0$ and  $r(\varepsilon) =\varepsilon^2$ if $\kappa=0$.

\item[(2)] For the homogenization regime (H2),
\begin{equation}\label{beta_e_H2}
\beta^e(k,\kappa, d, \varepsilon)= \dfrac{1}{\pi} \left(\ln \varepsilon + \ln 2 + \ln\dfrac{\pi}{d}\right) + \left(\dfrac{1}{2\pi}\sum_{n\neq 0} \dfrac{1}{|n|} - \dfrac{i}{d} \sum_{n=-\infty}^{\infty}  \dfrac{1}{\zeta_n(k) }\right),
\end{equation}
and
\begin{equation}\label{rho_e_H2}
\rho^e (X,Y; k,\kappa) = \dfrac{1}{\pi}  \ln (|X-Y|).
\end{equation}
In addition, $r(\varepsilon) = \varepsilon$ if $\kappa\neq 0$ and  $r(\varepsilon) = \varepsilon^2\ln\varepsilon$ if $\kappa=0$.
\end{enumerate}

\end{lem}

\medskip
\noindent\textbf{Remark 2.1}
\textit{ In the above, the subtraction
 $$\displaystyle{\frac{1}{2\pi}\sum_{n\neq 0} \frac{1}{|n|} - \frac{i}{d} \sum_{n=-\infty}^{\infty}\frac{1}{\zeta_n(k)  }}$$  is viewed as the sum of the converging series 
 $$\displaystyle{\sum_{n\neq 0} \left(\frac{1}{2\pi}\frac{1}{|n|} - \frac{i}{d}\frac{1}{\zeta_n(k)}\right) - \frac{i}{d}\frac{1}{\zeta_0(k) }}. $$ 
 Hence, the scalar function $\beta_e(k,\kappa, d, \varepsilon)$ is well defined.}  \\

\noindent\textbf{Proof}  
We derive the asymptotic expansion for the kernel $G_\varepsilon^e$ when $\kappa=0$.
For the homogenization regime (H1), we see that $kd\ll 1$. Therefore,
\begin{eqnarray*}
\sum_{n\neq 0} \dfrac{1}{\zeta_n(k)  } e^{i \kappa_n\varepsilon(X-Y)} &=& -\dfrac{id}{2\pi} \sum_{n\neq 0}\dfrac{1}{|n|\sqrt{1-(kd/2\pi n)^2}}e^{i \frac{2\pi n}{d}\varepsilon(X-Y)}  \\
 &=& -\dfrac{id}{2\pi} \sum_{n\neq 0}\dfrac{1}{|n|} \left( 1+\sum_{m=1}^\infty \dfrac{1\cdot3\cdots (2m-1)}{2^m m!} \left(\dfrac{kd}{2\pi n}\right)^{2m} \right) e^{i \frac{2\pi n}{d}\varepsilon(X-Y)} \\
&=&  -\dfrac{id}{2\pi} \sum_{n\neq0}\dfrac{1}{|n|} e^{i \frac{2\pi n}{d}\varepsilon(X-Y)}-\dfrac{id}{2\pi}  \sum_{m=1}^\infty \dfrac{1\cdot3\cdots (2m-1)}{2^m m!} \sum_{n\neq0} \left(\dfrac{kd}{2\pi n}\right)^{2m}\dfrac{1}{|n|} e^{i \frac{2\pi n}{d}\varepsilon(X-Y)} \\
&=& \dfrac{id}{2\pi}\ln\left(4\sin^2 (\pi \eta(X-Y))\right) +  O(k^2\varepsilon^3),
\end{eqnarray*}
where we have used the formula (cf. \cite{kress})
$$  \sum_{n\neq0}\dfrac{1}{|n|} e^{i \frac{2\pi n}{d}\varepsilon(X-Y)} = \ln\left(4\sin^2\dfrac{\pi \varepsilon(X-Y)}{d}\right). $$
Therefore, 
\begin{eqnarray*}
G_\varepsilon^e(X,Y) &=& -\dfrac{i}{d}\sum_{n=-\infty}^{\infty} \dfrac{1}{\zeta_n(k)  } e^{i \kappa_n\varepsilon(X-Y)} \\
&=& -\dfrac{i}{\zeta_0(k) d} + \dfrac{1}{2\pi}\ln\left(4\sin^2\dfrac{\pi \varepsilon(X-Y)}{d}\right) + O(k^2\varepsilon^2).
\end{eqnarray*}
The desired asymptotic expansion follows. 

For the homogenization regime (H2),  we have $k\ll1$ and $\varepsilon \ll 1$. Applying the Taylor expansion yields
\begin{eqnarray*}
\sum_{n\neq 0} \dfrac{1}{\zeta_n(k)   } e^{i \kappa_n\varepsilon(X-Y)} &=& -\dfrac{id}{2\pi} \sum_{n\neq 0}\dfrac{1}{|n|\sqrt{1-(kd/2\pi n)^2}}e^{i \frac{2\pi n}{d}\varepsilon(X-Y)}  \\
 &=& -\dfrac{id}{2\pi} \sum_{n\neq 0}\dfrac{1}{|n|} \left( 1+\sum_{m=1}^\infty \dfrac{1\cdot3\cdots (2m-1)}{2^m m!} \left(\dfrac{kd}{2\pi n}\right)^{2m} \right) e^{i \frac{2\pi n}{d}\varepsilon(X-Y)}.
\end{eqnarray*}
By the formula
$$ -\sum_{n\neq0}\dfrac{1}{|n|} e^{i \frac{2\pi n}{d}\varepsilon(X-Y)}= \ln\left(4\sin^2\dfrac{\pi \varepsilon(X-Y)}{d}\right), $$
and noting that for $m\ge1$ (cf. \cite{lewin}),
\begin{equation*}
 \sum_{n\neq 0}\dfrac{1}{|n|^{2m+1}}e^{i \frac{2\pi n}{d}\varepsilon(X-Y)}  = \sum_{n\neq 0 }\dfrac{1}{|n|^{2m+1}}+ O(\varepsilon^{2m}(X-Y)^{2m})\ln (\varepsilon(X-Y)),
\end{equation*}
we obtain 
\begin{eqnarray*}
\sum_{n\neq 0} \dfrac{1}{\zeta_n(k)} e^{i \kappa_n\varepsilon(X-Y)}  &=& \dfrac{id}{2\pi}\ln\left(4\sin^2\dfrac{\pi \varepsilon(X-Y)}{d}\right) \\
&& -\dfrac{id}{2\pi}  \sum_{m=1}^\infty \dfrac{1\cdot3\cdots (2m-1)}{2^m m!} \sum_{n\neq 0}  \left(\dfrac{kd}{2\pi}\right)^{2m}\dfrac{1}{|n|^{2m+1}} + O(\varepsilon^2\ln\varepsilon). \\
&=& \dfrac{id}{2\pi}\ln\left(4\sin^2\dfrac{\pi \varepsilon(X-Y)}{d}\right) -\dfrac{id}{2\pi}  \left( \sum_{n\neq 0} \dfrac{1}{|n|\sqrt{1-(kd/2\pi n)^2}}-\dfrac{1}{|n|} \right)+
O(\varepsilon^2\ln\varepsilon) \\
&=& \dfrac{id}{2\pi}\ln\left(4\sin^2\dfrac{\pi \varepsilon(X-Y)}{d}\right) + \dfrac{id}{2\pi} \sum_{n\neq 0}\dfrac{1}{|n|}  + \sum_{n\neq 0} \dfrac{1}{\zeta_n(k)}
+O(\varepsilon^2\ln\varepsilon).
\end{eqnarray*}
The desired asymptotic expansion follows by noting that $\varepsilon \ll d$ and using the expansion
$$\displaystyle{G_\varepsilon^e(X,Y) =-\dfrac{i}{d} \sum_{n=-\infty}^{\infty} \dfrac{1}{\zeta_n(k) } e^{i \kappa_n\varepsilon(X-Y)} = -\dfrac{i}{d} \left(\dfrac{1}{\zeta_0(k)} + \sum_{n\neq0} \dfrac{1}{\zeta_n(k)} e^{i \kappa_n\varepsilon(X-Y)}\right).}$$

Following a similar procedure, the asymptotic expansion of the kernel for $\kappa\neq0$ can be obtained, by noting that
 $$  -\dfrac{i}{d} \dfrac{1}{\zeta_0(k) } e^{i \kappa\varepsilon(X-Y)} =  -\dfrac{i}{\sqrt{k^2-\kappa^2} d} \left(1+i\kappa\varepsilon (X-Y)+\kappa^2 \cdot O(\varepsilon^2)\right). $$

 \qed

Let
\begin{equation}\label{beta}
\beta=\beta^i+\beta^e, 
\end{equation}
where $\beta^i$ is defined in \eqref{beta_i}, and $\beta^e$ is defined by \eqref{beta_e_H1} and \eqref{beta_e_H2} 
for two homogenization regimes respectively.
Set  
\begin{eqnarray*}
\rho(X,Y; k,\kappa) &=&\rho^i(X,Y)+\rho^e(X,Y; k,\kappa), \\
\rho_{\infty}(X,Y; k,\kappa)&=& r^i_{\varepsilon}(X,Y)+ r^e_{\varepsilon}(X,Y; k,\kappa),   \\
\tilde\rho_{\infty}(X,Y)&=& \tilde r_{\varepsilon}(X,Y).
\end{eqnarray*}
where $\rho^i$ is given by \eqref{rho_i},  $\rho^e$ is given by  \eqref{rho_e_H1} and \eqref{rho_e_H2}  for two homogenization regimes respecitvely, and $r^i_{\varepsilon}$, $\tilde r_{\varepsilon}$  and $r^e_{\varepsilon}$ are the high-order terms as specified in Lemma \ref{lem-green_slit} and Lemma \ref{lem-periodic_green}.
We define the integral operators $K$, $K_{\infty}$, $\tilde K_{\infty}$ by letting
\begin{eqnarray}\label{op_Ks}
  && (K \varphi)(X) = \int_0^1 \rho(X,Y; k,\kappa)  \varphi(Y) dY  \quad  X\in (0,1); \label{op_K} \\
 && (K_\infty\varphi) (X) = \int_0^1  \rho_{\infty} (X,Y; k,\kappa) \varphi(Y) dY  \quad  X\in (0,1); \label{op_K_inf} \\
 && (\tilde K_\infty\  \varphi)(X)  = \int_0^1 \tilde\rho_{\infty}(X, Y) \varphi(Y) dY  \quad  X\in (0,1).  \label{op_tK_inf}
\end{eqnarray}

\noindent\textbf{Remark 2.1} \textit{ Note that in the above, the function $\beta$,
the kernels of the integral operators $\rho$, $\rho_{\infty}$ and $\tilde\rho_{\infty}$  take different forms in the two homogenization regimes (H1) and (H2). Here and henceforth, we adopt the same notations for the sake of presenting a unified asymptotic framework for the scattering problem in two homogenization regimes (see Section \ref{sec-asym_scat_sol}). However, their values should be clear from the context. } \\

\begin{lem} \label{lem-optK}
Let the assumption in Lemma \ref{lem-periodic_green} hold, then in both homogenization regimes,
the operator $K$ is bounded from $V_1$ to $V_2$ with a bounded inverse.
Moreover, 
$$ \alpha(k,\kappa):=\langle K^{-1} 1, 1\rangle \; \mbox{is a real number and} \; \alpha(k,\kappa)\neq 0. $$
\end{lem}

\medskip

\noindent\textbf{Remark 2.2}
$\alpha$ takes different values in the two homogenization regimes (H1) and (H2).
It depends on $k, \kappa$ in the former homogenization regime, and is independent of $k$ and $\kappa$ in the latter.
For the ease of notation, we will simply denote it as $\alpha$ in the rest of the paper. \\

\medskip

\noindent\textbf{Proof} The proof for the homogenization regime (H1) is postponed to  the appendix. 
\medskip
For the homogenization regime (H2), recall the kernel of the $K$ takes the form
$$ \rho(X,Y; k,\kappa) =\dfrac{1}{\pi}  \ln |X-Y|  + \dfrac{1}{\pi} \left[ \ln \left(\abs{\sin \left(\frac{\pi(X+Y)}{2}\right)}\right) + \ln \left(\abs{\sin \left(\frac{\pi(X-Y)}{2}\right)}\right) \right].  $$
Note that $ \rho$ is independent of $k, \kappa$.
The proof can be found in Theorem 4.1 and Lemma 4.2 of \cite{eric10}.
\medskip

\begin{lem} \label{lem-operators}
Let the assumption in Lemma  \ref{lem-periodic_green} hold, then the following holds for (H1) and (H2).
\begin{enumerate}
\item[(1)]
The operator $T^e+T^i $ admits the following decomposition:
$$ T^e+T^i = \beta P + K + K_\infty. $$
Moreover, $K_\infty$ is bounded from $V_1$ to $V_2$ with the operator norm $\| K_{\infty}\|  \lesssim r(\varepsilon)$ uniformly for bounded $k$'s.

\item[(2)]
The operator  $\tilde T^i$ admits the following decomposition: 
$$\tilde T^i = \tilde \beta P + \tilde K_\infty,$$ 
Moreover, $\tilde K_\infty$ is bounded from $V_1$ to $V_2$ with the operator norm
$\| K_{\infty}\|  \lesssim e^{-1/\varepsilon}$ uniformly for bounded $k$'s. 

\end{enumerate}
\end{lem}

\noindent\textbf{Proof}   First, from the definition of the operators $T^e$ and $T^i$, the kernel of $T^e+T^i $ 
is $G_\varepsilon^e(X, Y)+G_\varepsilon^i(X, Y)$. By the asymptotic expansions of the kernels in Lemma \ref{lem-green_slit} and \ref{lem-periodic_green}, it is clear that
\begin{eqnarray*}
G_\varepsilon^e(X, Y)+G_\varepsilon^i(X, Y) &=& \beta^e+\beta^i + \rho^e(X,Y;k,\kappa) + \rho^i(X,Y) + r^e_\varepsilon(X,Y;k,\kappa) +  r^i_\varepsilon(X,Y) \\
&=& \beta + \rho(X,Y;k,\kappa)  + \rho_{\infty}(X,Y).
\end{eqnarray*}
The assertion (1) follows by the definition of the operator \eqref{op_K} and \eqref{op_K_inf}. The proof of (2) follows by using the decomposition
$$ \tilde G_\varepsilon^i(X, Y) =  \tilde\beta(k,\varepsilon) + \tilde r_{\varepsilon}((X, Y) = \tilde \beta(k,\varepsilon) + \tilde \rho_{\infty}(X, Y). $$

\subsection{Asymptotic expansion of the solution to the scattering problem}\label{sec-asym_scat_sol}

For both homogenization regimes, we define
\begin{equation*}
\mathbb{P}= \left[
\begin{array}{cc}
\beta P    & \tilde \beta P \\
\tilde \beta P  & \beta P 
\end{array}
\right],
\quad
\mathbb{K}_\infty=
\left[
\begin{array}{cc}
K_\infty    & \tilde K_\infty \\
\tilde K_\infty  & K_\infty
\end{array}
\right], 
\quad
\mathbf{f}=\left[ \begin{array}{c}
f/\varepsilon \\
0
\end{array}\right],
\quad\mbox{and}\quad
\mathbb{L}=K \mathbb{I}  +\mathbb{K}_\infty.
\end{equation*}
Then from the decomposition of the operators in Lemma \ref{lem-operators}, we may rewrite the system of the integral equations
\eqref{eq-scattering3} as
\begin{equation}\label{eq-scattering4}
(\mathbb{P} + \mathbb{L}) \boldsymbol{\varphi} =   \mathbf{f}. 
\end{equation}
Next, we derive the asymptotic expansion of the solution $\boldsymbol{\varphi}$.
By Lemma \ref{lem-optK}, it is also easy to see that $\mathbb{L}$ is invertible for sufficiently small 
$\varepsilon$. 
Applying the Neumann series yields
$$ \mathbb{L}^{-1}  = \left(K\mathbb{I}+\mathbb{K}_\infty\right)^{-1} = \left(\sum_{j=0}^\infty (-1)^j\left(K^{-1}\mathbb{K_\infty} \right)^j \right)  K^{-1} = K^{-1}\mathbb{I} + O\left(r(\varepsilon)\right). $$
Therefore, the following lemma follows immediately.
\begin{lem} \label{lem-L_inv} 
 Let $\mathbf{e}_1 = [1, 0]^T $ and $\mathbf{e}_2 = [0, 1]^T $. Then
\begin{equation} \label{eq-II}
\mathbb{L}^{-1}  \mathbf{e}_1 = K^{-1}1 \cdot \mathbf{e}_1 + O(r(\varepsilon)), \quad
\mathbb{L}^{-1}  \mathbf{e}_2 = K^{-1}1 \cdot \mathbf{e}_2 + O(r(\varepsilon)),
\end{equation}
and 
\begin{equation} \label{eq-III}
\langle \mathbb{L}^{-1}  \mathbf{e}_1,  \mathbf{e}_1 \rangle = \alpha +  O(r(\varepsilon)), \quad 
\langle \mathbb{L}^{-1}  \mathbf{e}_1,  \mathbf{e}_2\rangle =O(r(\varepsilon)).
\end{equation}
Here $\alpha$ is defined in Lemma \ref{lem-optK}.
\end{lem}

\begin{lem} \label{lem-identity1}
 Let $\mathbf{e}_1 = [1, 0]^T $ and $\mathbf{e}_2 = [0, 1]^T $. Then
\begin{equation*} \label{eq-l} \langle \mathbb{L}^{-1}  \mathbf{e}_1,  \mathbf{e}_1 \rangle = \langle \mathbb{L}^{-1}  \mathbf{e}_2,  \mathbf{e}_2 \rangle, \quad \langle \mathbb{L}^{-1}  \mathbf{e}_1,  \mathbf{e}_2 \rangle = \langle   \mathbb{L}^{-1}\mathbf{e}_2,  \mathbf{e}_1 \rangle.
\end{equation*}
\end{lem}

\noindent\textbf{Proof} Let $\mathbb{L}^{-1}  \mathbf{e}_1 = (a, b)^T$. Then $\mathbb{L}(a, b)^T =\mathbf{e}_1 $. More precisely, 
\begin{eqnarray*}
K a + K_{\infty} a + \tilde K_\infty b &=&1, \\
K b + \tilde K_{\infty} a +  K_\infty b &=&0.
\end{eqnarray*}
It follows that $\mathbb{L}(b, a)^T =\mathbf{e}_2 $, or equivalently, 
$$
\mathbb{L}^{-1}  \mathbf{e}_2 = (b, a)^T,
$$
hence the two identities hold. \qed \\

\medskip

By applying $\mathbb{L}^{-1}$ on both sides of \eqref{eq-scattering4},  we see that
\begin{equation}\label{eq-scattering5}
\mathbb{L}^{-1} \;\mathbb{P} \; \boldsymbol{\varphi} + \boldsymbol{\varphi} =  \mathbb{L}^{-1}  \mathbf{f}.
\end{equation}
Note that
 $$  \mathbb{P} \; \boldsymbol{\varphi} = \beta  \langle  \boldsymbol{\varphi}, \mathbf{e}_1 \rangle \mathbf{e}_1  + \beta \langle  \boldsymbol{\varphi}, \mathbf{e}_2 \rangle \mathbf{e}_2
 + \tilde \beta \langle  \boldsymbol{\varphi}, \mathbf{e}_2 \rangle \mathbf{e}_1  + \tilde \beta \langle  \boldsymbol{\varphi}, \mathbf{e}_1 \rangle \mathbf{e}_2,$$
the above operator equation can be written as
\begin{equation}\label{Op_eqns_1}
\beta \langle  \boldsymbol{\varphi}, \mathbf{e}_1 \rangle  \mathbb{L}^{-1}  \mathbf{e}_1  + \beta  \langle  \boldsymbol{\varphi}, \mathbf{e}_2 \rangle  \mathbb{L}^{-1}  \mathbf{e}_2 + 
\tilde \beta \langle  \boldsymbol{\varphi}, \mathbf{e}_2 \rangle \mathbb{L}^{-1} \mathbf{e}_1  + \tilde \beta \langle  \boldsymbol{\varphi}, \mathbf{e}_1 \rangle \mathbb{L}^{-1} \mathbf{e}_2+
\boldsymbol{\varphi} =   \mathbb{L}^{-1}  \mathbf{f}. 
\end{equation}
By taking the inner product of \eqref{Op_eqns_1} with $\mathbf{e}_1$ and $\mathbf{e}_2$ respectively, it follows that
\begin{equation}\label{eq-linear_sys}
( \mathbb{M}+\mathbb{I})\left[
 \begin{array}{llll}
\langle  \boldsymbol{\varphi}, \mathbf{e}_1 \rangle  \\
\langle  \boldsymbol{\varphi}, \mathbf{e}_2 \rangle
\end{array}
\right] =
\left[
\begin{array}{llll}
 \langle \mathbb{L}^{-1} \mathbf{f}, \mathbf{e}_1 \rangle   \\
 \langle \mathbb{L}^{-1} \mathbf{f}, \mathbf{e}_2 \rangle 
\end{array}
\right],
\end{equation}
where the matrix $\mathbb{M}$ is defined as
\begin{equation}  \label{eq-matrix-m}
\mathbb{M}:=
\beta  \left[
\begin{array}{llll}
 \langle \mathbb{L}^{-1}  \mathbf{e}_1,  \mathbf{e}_1 \rangle  &  \langle \mathbb{L}^{-1}  \mathbf{e}_2,  \mathbf{e}_1\rangle \\
 \langle \mathbb{L}^{-1}  \mathbf{e}_1,  \mathbf{e}_2 \rangle  &  \langle \mathbb{L}^{-1}  \mathbf{e}_2,  \mathbf{e}_2 \rangle 
 \end{array}
\right] 
+
\tilde \beta \left[
\begin{array}{llll}
 \langle \mathbb{L}^{-1}  \mathbf{e}_2,  \mathbf{e}_1 \rangle  &  \langle \mathbb{L}^{-1}  \mathbf{e}_1,  \mathbf{e}_1\rangle \\
 \langle \mathbb{L}^{-1}  \mathbf{e}_2,  \mathbf{e}_2 \rangle  &  \langle \mathbb{L}^{-1}  \mathbf{e}_1,  \mathbf{e}_2 \rangle 
 \end{array}
\right]. 
\end{equation}
From Lemma \ref{lem-identity1}, it is observed that
\begin{equation*}
\mathbb{M} = 
\left(\beta  +
\tilde \beta 
\left[
\begin{array}{llll}
0 & 1 \\
1 & 0 
\end{array}
\right]\right)
\left[
\begin{array}{llll}
 \langle \mathbb{L}^{-1}  \mathbf{e}_1,  \mathbf{e}_1 \rangle  &  \langle \mathbb{L}^{-1}  \mathbf{e}_1,  \mathbf{e}_2\rangle \\
 \langle \mathbb{L}^{-1}  \mathbf{e}_1,  \mathbf{e}_2 \rangle  &  \langle \mathbb{L}^{-1}  \mathbf{e}_1,  \mathbf{e}_1 \rangle 
 \end{array}
\right].
\end{equation*}
A straightforward calculation shows that the eigenvalues of $\mathbb{M}+\mathbb{I}$ are
\begin{eqnarray}\label{eq-eigen_M}
\lambda_1(k; \kappa, d, \varepsilon) &=& 1+(\beta+ \tilde \beta )   \left(\langle \mathbb{L}^{-1}  \mathbf{e}_1,  \mathbf{e}_1 \rangle  + \langle \mathbb{L}^{-1}  \mathbf{e}_1,  \mathbf{e}_2\rangle\right), \label
{eq-lambda1}\\
\lambda_2(k; \kappa, d, \varepsilon) &=& 1+(\beta-\tilde \beta )   \left(\langle \mathbb{L}^{-1}  \mathbf{e}_1,  \mathbf{e}_1 \rangle  - \langle \mathbb{L}^{-1}  \mathbf{e}_1,  \mathbf{e}_2\rangle\right),
\label{eq-lambda2}
\end{eqnarray}
and the associated eigenvectors are $[1 \quad 1]^T$ and $[1 \quad -1]^T$.  
For simplicity of notation, let us define 
\begin{equation}\label{eq-def_pq}
p(k; \kappa, d, \varepsilon):=\varepsilon \lambda_1(k; \kappa, d, \varepsilon) \quad \mbox{and} \quad q(k; \kappa, d, \varepsilon):=\varepsilon \lambda_2(k; \kappa, d, \varepsilon),
\end{equation}
which will be used throughout the rest of the paper.

Solving \eqref{eq-linear_sys} leads to
\begin{equation}\label{eq-phi_dot_e1e2}
\left[
\begin{array}{llll}
\langle  \boldsymbol{\varphi}, \mathbf{e}_1 \rangle  \\
\langle  \boldsymbol{\varphi}, \mathbf{e}_2 \rangle
\end{array}
\right] =( \mathbb{M}+\mathbb{I})^{-1} 
\left[
\begin{array}{llll}
 \langle \mathbb{L}^{-1} \mathbf{f}, \mathbf{e}_1 \rangle   \\
 \langle \mathbb{L}^{-1} \mathbf{f}, \mathbf{e}_2 \rangle 
\end{array}
\right].
\end{equation}
By substituting into \eqref{eq-scattering5}, we obtain the solution of the integral equation system \eqref{eq-scattering3}:
\begin{equation}\label{eq-phi}
\boldsymbol{\varphi} = \mathbb{L}^{-1}  \mathbf{f}
- \bigg[  \mathbb{L}^{-1}  \mathbf{e}_1 \quad  \mathbb{L}^{-1}  \mathbf{e}_2 \bigg]  
\left[
\begin{array}{llll}
\beta  &   \tilde \beta \\
\tilde \beta  &\beta
\end{array}
\right] 
(\mathbb{M}+\mathbb{I})^{-1} 
\left[
\begin{array}{llll}
 \langle \mathbb{L}^{-1} \mathbf{f}, \mathbf{e}_1 \rangle   \\
 \langle \mathbb{L}^{-1} \mathbf{f}, \mathbf{e}_2 \rangle 
\end{array}
\right].
\end{equation}

\begin{lem} \label{lem-phi}
Assume that $k\in\mathbf{R}^+$ is not an eigenvalue of the scattering operator. 
Let $\kappa= k\sin \, \theta$, where $\theta$ is the incident angle.
Then the following asymptotic expansion holds for the solution $\boldsymbol{\varphi}$ of \eqref{eq-scattering3} in $V_1 \times V_1$
in both homogenization regimes:
\begin{eqnarray} \label{eq-varphi}
 \boldsymbol{\varphi} &=& K^{-1}1 
 \cdot  \left[ \kappa \cdot O(1) \cdot\mathbf{e}_1 + \dfrac{\alpha}{p}  (\mathbf{e}_1 + \mathbf{e}_2 ) + \dfrac{\alpha}{q} (\mathbf{e}_1 - \mathbf{e}_2)  \right] \nonumber + \left(\dfrac{\alpha}{p} +  \dfrac{\alpha}{q} \right) \cdot O(r(\varepsilon)) +  O(r(\varepsilon)).
\end{eqnarray}
Moreover, 
\begin{equation}  \label{eq-varphi-1}
\left[
\begin{array}{llll}
\langle  \boldsymbol{\varphi}, \mathbf{e}_1 \rangle  \\
\langle  \boldsymbol{\varphi}, \mathbf{e}_2 \rangle
\end{array}
\right]  
= \bigg[\alpha+O(r(\varepsilon)) \bigg] \left(
\dfrac{1}{p}
\left[
\begin{array}{cc}
1  \\
1 
\end{array}
\right] 
+
\dfrac{1}{q}
\left[
\begin{array}{cc}
1 \\
-1 
\end{array}
\right] \right).
\end{equation}
Here $\alpha$ is defined in Lemma \ref{lem-optK}, and $p$ and $q$ are defined by \eqref{eq-def_pq}.
\end{lem}

\noindent\textbf{Proof}.  For given $k$ and $\kappa$, we see that $\zeta:=\sqrt{k^2-\kappa^2} = k\cos\theta$. Thus 
 $\kappa/\zeta  = \tan\theta$ is bounded  and the assumption in Lemma \ref{lem-periodic_green} and \ref{lem-operators} holds. 
By applying the asymptotic expansion derived in the previous section, we obtain the representations
 \eqref{eq-phi_dot_e1e2} and \eqref{eq-phi} for $\boldsymbol{\varphi}$, $\langle  \boldsymbol{\varphi}, \mathbf{e}_1 \rangle$,
 $\langle  \boldsymbol{\varphi}, \mathbf{e}_2 \rangle$.
 
Note that the matrix $\mathbb{M}+\mathbb{I}$ has two eigenvalues $\lambda_1$ and $\lambda_2$ given by \eqref{eq-lambda1} and \eqref{eq-lambda2}, 
which are associated with the eigenvectors $[1 \quad 1]^T$ and $[1 \quad -1]^T$ respectively, it follows that
\begin{equation*}
(\mathbb{M}+\mathbb{I})^{-1}  =
\dfrac{1}{2\lambda_1}
\left[
\begin{array}{cc}
1 & 1 \\
1 & 1 
\end{array}
\right] 
+
\dfrac{1}{2\lambda_2}
\left[
\begin{array}{cc}
1 & -1 \\
-1 & 1 
\end{array}
\right].
\end{equation*}
Substituting into \eqref{eq-phi_dot_e1e2} and \eqref{eq-phi} yields
\begin{eqnarray*}
\left[
\begin{array}{llll}
\langle  \boldsymbol{\varphi}, \mathbf{e}_1 \rangle  \\
\langle  \boldsymbol{\varphi}, \mathbf{e}_2 \rangle
\end{array}
\right] 
&=&
\dfrac{1}{2\lambda_1}\langle \mathbb{L}^{-1} \mathbf{f}, \mathbf{e}_1+\mathbf{e}_2 \rangle
\left[
\begin{array}{cc}
1  \\
1 
\end{array}
\right] 
+
\dfrac{1}{2\lambda_2(k,\varepsilon)}\langle \mathbb{L}^{-1} \mathbf{f}, \mathbf{e}_1-\mathbf{e}_2 \rangle
\left[
\begin{array}{cc}
1 \\
-1 
\end{array}
\right],
\end{eqnarray*}
and 
\begin{eqnarray}\label{eq-varphi-2}
\boldsymbol{\varphi} &=& \mathbb{L}^{-1}  \mathbf{f}
+ \dfrac{1- \lambda_1/\langle\mathbb{L}^{-1}\mathbf{e}_1, \mathbf{e}_1+ \mathbf{e}_2 \rangle}{2\lambda_1} 
\langle \mathbb{L}^{-1} \mathbf{f}, \mathbf{e}_1 + \mathbf{e}_2  \rangle \cdot 
(\mathbb{L}^{-1}  \mathbf{e}_1 + \mathbb{L}^{-1}  \mathbf{e}_2) \nonumber  \\
&& + \dfrac{1- \lambda_2(k,\varepsilon)/\langle\mathbb{L}^{-1}\mathbf{e}_1, \mathbf{e}_1+ \mathbf{e}_2 \rangle}{2\lambda_2(k,\varepsilon)} 
\langle \mathbb{L}^{-1} \mathbf{f}, \mathbf{e}_1 - \mathbf{e}_2 \rangle \cdot
(\mathbb{L}^{-1}  \mathbf{e}_1 - \mathbb{L}^{-1}  \mathbf{e}_2).
\end{eqnarray}

From the Taylor expansion of  $\mathbf{f}$ and the asymptotic expansion of the operator $\mathbb{L}^{-1}$ in Lemma \ref{lem-L_inv}, we can obtain the desired asymptotic expansions for $\boldsymbol{\varphi}$, $\langle  \boldsymbol{\varphi}, \mathbf{e}_1 \rangle$,
 $\langle  \boldsymbol{\varphi}, \mathbf{e}_2 \rangle$. This derivation is similar to Lemma 3.5 of the first part of this series \cite{lin_zhang17},
 and we omit here. \qed \\

\subsection{An overview of diffraction anomaly and field enhancement}
From Lemma \ref{lem-phi}, we see that the solution of the system of integral equations $\boldsymbol{\varphi}$ depends on the values of
two functions $p$ and $q$. In the rest of the paper, we investigate their values in two homogenization regimes, and explore anomalous behaviors
and field enhancement for the solution to the scattering problem.

In the homogenization regime (H1), we will shown that for each $\kappa$, there exists roots for $p(k; \kappa, d, \varepsilon)=0$ and $q(k; \kappa, d, \varepsilon)=0$ such that the homogeneous scattering problem attains nontrivial solutions. Indeed, such roots correspond to
the eigenvalues of the scattering operator, and very interestingly, the first branch of the dispersion curve $k(\kappa)$,
and the corresponding localized eigenmodes resemble those of surface plasmon polaritons of the nobel metal slab.
This is so-called spoof surface plasmon effect, which mimics surface plasmon of noble metals in a perfect conductor by corrugating its surface \cite{garcia05, pendry04}. It extends the frequency band for the surface plasmon, which is originally supported on a flat noble metal in optical and near-infrared regime, to the 
terahertz or lower frequency regime where metals are nearly perfect conductors.
We will derive the asymptotic expansions for the dispersion curve and the associated eigenmodes in Section \ref{sec-H1}.
A discussion of the surface plasmonic effect will also be presented.

It is also demonstrated, to our surprise, that total transmission can be achieved in this homogenization regime. More precisely, for an incident plane wave, there exist certain frequencies such that no wave is reflected, and all electromagnetic energy passes through the slab in the limiting case of $\varepsilon\to0$. Such phenomenon also occurs for all frequencies at a specific incident angle. These results will also be reported in Section \ref{sec-H1}.

In the homogenization regime (H2), it will be shown that although no roots exist for $p(k; \kappa, d, \varepsilon)=0$ and $q(k; \kappa, d, \varepsilon)=0$,
the values of $p$ and $q$ will contribute to in a way that leads to significant electric field enhancement in and near the slits. The 
asymptotic expansions of the electromagnetic fields will be derived and their enhancement behaviors will be investigated in Section \ref{sec-H2}.

\section{Homogenization regime (H1): surface bound-state modes and total transmission}\label{sec-H1}
In the homogenization regime (H1),  the scaling of parameters are given by
$\varepsilon \sim d \ll 1$ (Figure \ref{fig-geo_H1H2}, top).  
It is known a nobel metal slab support surface plasmonic waves in the optical and near-infrared regime, 
but such localized plasmonic waves do not exist at lower frequencies when the metal is close to a perfect conductor \cite{maier07}. 
In this section, we demonstrate that when the perfect conducting slab
is perforated by an array of small slits and with small period as shown in Figure \ref{fig-geo_H1H2} (top), 
then the associated dispersion curve would resemble that of surface plasmon 
polaritons of the nobel metal slab. In addition, surface bound states, which resemble the plasmonic waves, are supported on top of the perfect conducting slab. To this end, we derive the asymptotic expansions of the dispersion relation and the corresponding eigenmodes in Section \ref{sec-H1_asympotics_dispersion} and \ref{sec-H1_eigenmode}. The effective medium in the slab induced by the periodic structure  as $\varepsilon \to 0$ is derived in Section \ref{sec-H1_effective_medium}, which recovers the leading order of the dispersion relation given in Section \ref{sec-H1_asympotics_dispersion}.
A brief discussion on the surface plasmon effect of the perfect conducting conducting slab with slits
and that of the plasmonic metal is given in Section \ref{sec-H1_spoof_spp}.

The other phenomenon induced by the given periodic structure is the total transmission through the small slits when an incident plane wave impinges on the slab. This occurs either at certain frequencies for all incident angles or all frequencies for a specific incident angle. More precisely, no wave is reflected, and all electromagnetic energy passes through the slab in the limiting case of $\varepsilon\to0$. We derive the field pattern above and below the slab for
the scattering problem in Section \ref{sec-H1_effective_medium} and discuss the total transmission phenomenon in Section \ref{sec-H1_total_trans}.

\subsection{Asymptotic expansions of the dispersion relation}\label{sec-H1_asympotics_dispersion}
To obtain the dispersion relation, we consider the homogeneous scattering problem wherein the incident wave $u^i=0$.
By \eqref{eq-scattering4}, the homogeneous problem is equivalent to the operator equation
$$ (\mathbb{P} + \mathbb{L})\boldsymbol{\varphi}=0. $$
In light of \eqref{eq-linear_sys}, this reduces to
\begin{equation*}
( \mathbb{M}+\mathbb{I})\left[
 \begin{array}{llll}
\langle  \boldsymbol{\varphi}, \mathbf{e}_1 \rangle  \\
\langle  \boldsymbol{\varphi}, \mathbf{e}_2 \rangle
\end{array}
\right] =0,
\end{equation*}
where the matrix $\mathbb{M}$ is defined by \eqref{eq-matrix-m}.
Therefore, the characteristic values of the operator-valued function $\mathbb{P} + \mathbb{L}$,
or equivalently the eigenvalues of the scattering operator,
 are the roots of $\lambda_1(k; \kappa, d, \varepsilon) $ and $\lambda_2(k; \kappa, d, \varepsilon)$, the eigenvalues of $\mathbb{M}+\mathbb{I}$.
 Then one only needs to solve $p(k; \kappa, d, \varepsilon) =0$ and $q(k; \kappa, d, \varepsilon) =0$ 
to obtain the eigenvalues of the scattering operator.

In light of \eqref{eq-lambda1}, and the definition of $\beta$ in \eqref{beta}, and Lemma \ref{lem-L_inv},
we may explicitly express $p$ as follows:
\begin{eqnarray}\label{eq-p_H1}
p(k; \kappa, d, \varepsilon)  &=& \varepsilon+\left[ \left(\dfrac{\cot k }{k} + \dfrac{1}{k\sin k} - \dfrac{i \eta}{\sqrt{k^2-\kappa^2}} \right)+\dfrac{3\ln 2}{\pi} \varepsilon  \right] 
  \left(\langle \mathbb{L}^{-1}  \mathbf{e}_1,  \mathbf{e}_1 \rangle  + \langle \mathbb{L}^{-1}  \mathbf{e}_1,  \mathbf{e}_2\rangle\right) \nonumber \\
&=& 
\varepsilon+\left[ \left(\dfrac{\cot k }{k} + \dfrac{1}{k\sin k} - \dfrac{i \eta}{\sqrt{k^2-\kappa^2}} \right)+\dfrac{3\ln 2}{\pi} \varepsilon  \right]  \left(\alpha + s(\varepsilon) \right),
\end{eqnarray}
where $s(\varepsilon)\sim O(r(\varepsilon))$. Similarly,
\begin{equation}\label{eq-q_H1}
q(k; \kappa, d, \varepsilon)=\varepsilon+\left[ \left(\dfrac{\cot k }{k} - \dfrac{1}{k\sin k} - \dfrac{i \eta}{\sqrt{k^2-\kappa^2}} \right)+\dfrac{3\ln 2}{\pi} \varepsilon \right] \left(\alpha + t(\varepsilon) \right).
\end{equation}
where $t(\varepsilon)\sim O(r(\varepsilon))$. First, we investigate the roots for the leading-order terms of $p$ and $q$.

\begin{lem}\label{lem-disp_curve}
For each $\kappa$,
$$ c_\pm(k,\kappa)=\dfrac{\cot k }{k} \pm \dfrac{1}{k\sin k} - \dfrac{i \eta}{\sqrt{k^2-\kappa^2}}=0 $$
attains real roots $k_{m,0}^\pm(\kappa)$ ($m=0, 1, 2, \cdots, M^{\pm}$). In addition, 
\begin{itemize}
\item[(i)] $ 0\le k_{0,0}^\pm(\kappa) < k_{1,0}^\pm(\kappa) < \cdots < k_{M^\pm,0}^\pm(\kappa) \le |\kappa|$. 
\item[(ii)] For each $m$, $k_{m,0}^\pm(\kappa)$ is a continuous and monotonic function of $\kappa$.
\item[(iii)]  As $|\kappa|\to\infty$, $k_{m,0}^+(\kappa) \to m\pi$ and $k_{m,0}^-(\kappa) \to (m+1)\pi$ if $m$ is odd,
and $k_{m,0}^+(\kappa) \to (m+1)\pi$ and  $k_{m,0}^-(\kappa) \to (m+2)\pi$ if $m$ is even.
\end{itemize}
\end{lem}

\noindent\textbf{Proof}. Solving $c_+(k,\kappa) =0$ yields
$$ \kappa = \pm\,\phi_+(k):=\pm \, k \sqrt{1+ \eta^2 \tan^2(k/2)}, \quad k\ge 0.$$
Without loss of generality, we consider $\kappa\ge0$ and  $\kappa =\phi_+(k)$.

\begin{figure}[!htbp]
\begin{center}
\includegraphics[height=4.5cm,width=8cm]{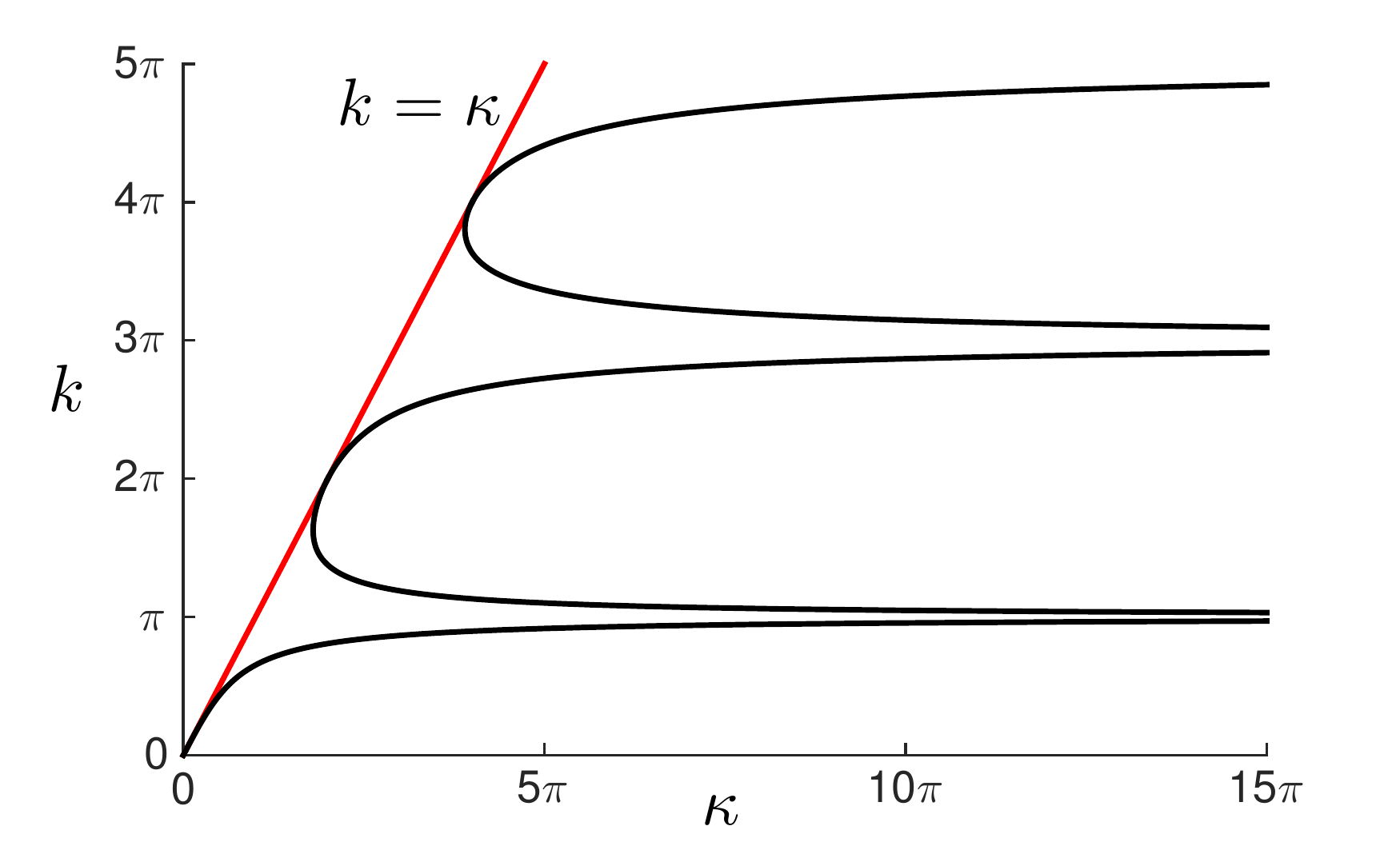}
\includegraphics[height=4.5cm,width=8cm]{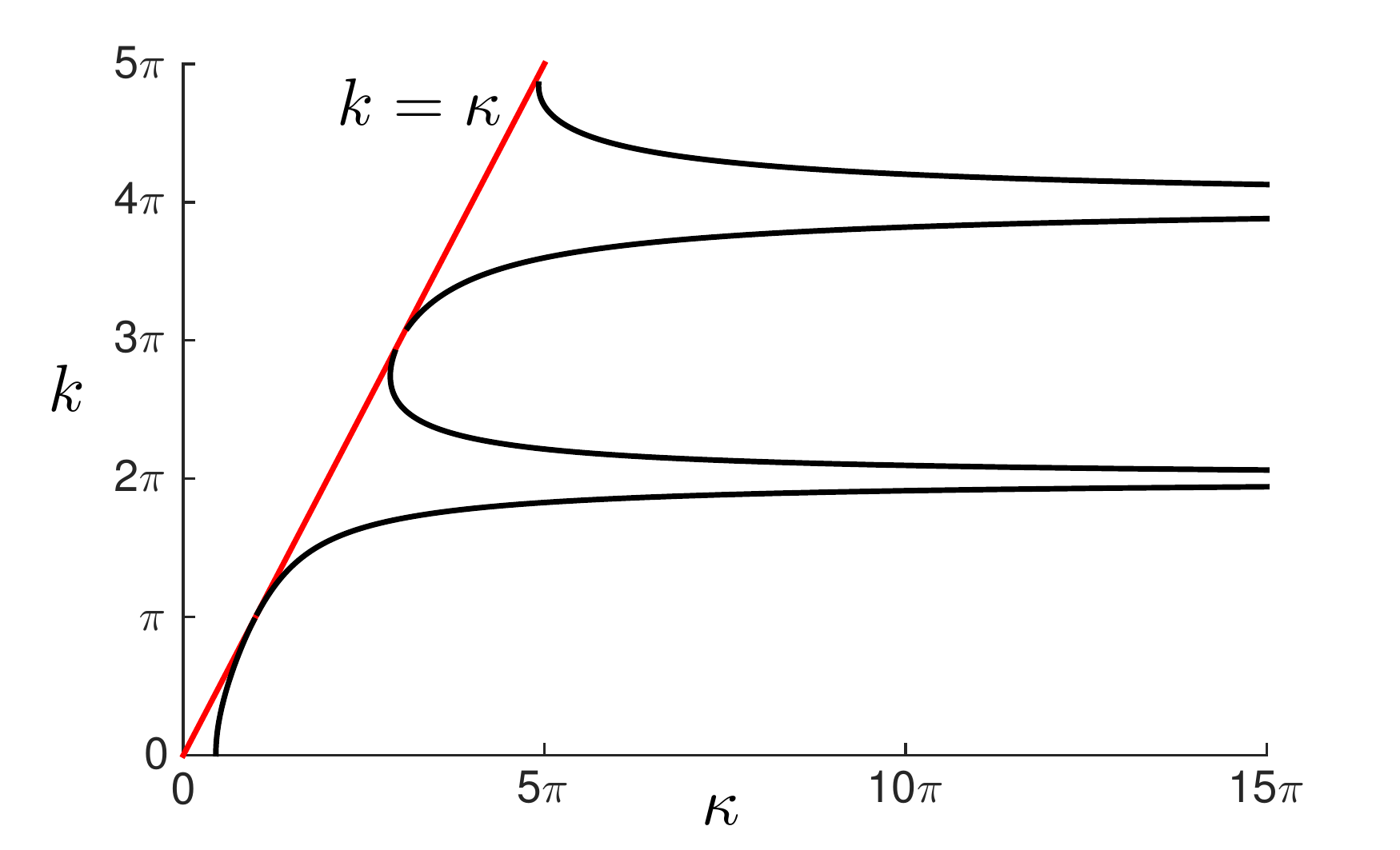}
\caption{Root of $c_+(k,\kappa) =0$: $\kappa=\phi^+(k)=k \sqrt{1+ \eta^2 \tan^2(k/2)}$ (left),
and root of $c_-(k,\kappa) =0$: $\kappa=\phi^-(k)=k \sqrt{1+ \eta^2 \cot^2(k/2)}$ (right).}\label{H1_dispersion}
\end{center}
\end{figure}

Decompose the domain of the definition for $\phi_+(k)$ as non-overlapping intervals:
$$D(\phi_+)=\displaystyle{\bigcup_{m=0}^{\infty} \Big(\big[2m\pi, (2m+1)\pi\big) \bigcup \big((2m+1)\pi, (2m+2)\pi\big)\Big)}.$$
Then for $k\in \big[2m\pi, (2m+1)\pi\big)$, $\phi_+(k)$ is a monotonic increasing and its range is $\big[2m\pi,+\infty\big)$ (cf. Figure \ref{H1_dispersion}).
Therefore, the inverse
$$(\phi_+)^{-1}: \quad \big[2m\pi,+\infty\big) \to \big[2m\pi, (2m+1)\pi\big) $$
exists, which we denote by $k_{2m,0}^+(\kappa)$.  
It is clear that $k_{2m,0}^+(\kappa)$ is continuous and monotonic. Furthermore,
$k_{2m,0}^+(\kappa)\in \big[2m\pi, (2m+1)\pi\big)$ and $k_{2m,0}^+(\kappa)\le \kappa$.
As $\kappa\to\infty$,  it follows that $k_{2m,0}^+(\kappa)\to (2m+1)\pi$.

Similarly, $\phi_+(k)$ is a monotonic decreasing in the interval 
$\big((2m+1)\pi, (2m+2)\pi\big)$ with range $\big((2m+2)\pi,+\infty\big)$ (cf. Figure \ref{H1_dispersion}). The inverse $(\phi_+)^{-1}$ also exists
and is denoted by $k_{2m+1,0}^+(\kappa)$.
We have $k_{2m+1,0}^+(\kappa)\in \big((2m+1)\pi, (2m+2)\pi\big)$ and $|k_{2m+1,0}^-(\kappa)| < \kappa$.
The continuity, monotonicity and the asymptotic behavior of the function are straighforward to derive.

Since the range of $k_{2m,0}^+(\kappa)$ and $k_{2m+1,0}^+(\kappa)$ does not overlap for different values of $m$, 
we may arrange the roots such that
$ 0\le k_{0,0}^+ < k_{1,0}^+ < \cdots < k_{M,0}^+\le \kappa. $
Similarly, by solving $c_-(k,\kappa) =0$ we obtain
$$ \kappa = \pm\,\phi_-(k):=\pm \, k \sqrt{1+ \eta^2 \cot^2(k/2)}, \quad k\ge 0.$$
An analogous argument as above leads to the assertion for the roots of $c_-(k,\kappa)$. \qed \\

Next we derive the asymptotic expansion of the roots for $p$ and $q$. Note that $kd\ll1$ in the homogenization regime (H1),
we may restrict the discussion in the bounded domain $D_M:=\{ z: \; |z| \le M\}$ on the complex $k$-plane, where $M>0$ is a fixed constant. 
In addition, for a given $\kappa$, $p$ and $q$ are analytic with respect to $k$ in $D_M$ except for the cut-off frequency $k=\kappa$, 
thus we consider $k$ away from such cut-off frequency. To this end,
 let us define the domain
$$ D_{\kappa,\delta, M} := \{ z: \; |z| \le M\} \backslash B_{\delta}(\kappa), $$
where $\delta$ is a positive constant and $B_{\delta}(z)$ is the disk with radius $\delta$ centered at $z$ on the complex plane.
Let $k_{m,0}^{\pm}$ be the roots of $c_\pm(k,\kappa)=0$ as given in Lemma \ref{lem-disp_curve}.
Note that $\partial_k c_{\pm}(k_{m,0}^{\pm},\kappa) = 0$ only if $\kappa=0$ or $\partial_k\phi^{\pm}(k_{m,0})=0$.
They hold on a countable set on the $(\kappa,k)$-plane, as observed from the definition of $\phi^{\pm}(k)$ and Figure \ref{H1_dispersion}.
If $k_{m,0}\in D_{\kappa,\delta, M}$, we obtain the following asymptotic expansion in the neighborhood of $k_{m,0}^{\pm}$.

\begin{thm}\label{thm-asym_dispersion}
For each $\kappa$, if $k_{m,0}\in D_{\kappa,\delta, M}$ and $\partial_k c_{\pm}(k_{m,0}^{\pm},\kappa)\neq 0$.
Then in the neighborhood  $B_{\delta/2}(k_m^{\pm})$,
the roots of $p(k; \kappa, d, \varepsilon) =0$ and $q(k; \kappa, d, \varepsilon) =0$
attain the following asymptotic expansion :
\begin{equation}\label{eq-asym_res_eig}
k_m^{\pm}= k_m^{\pm}(\kappa, \varepsilon)=  k_{m,0}^{\pm}(\kappa)+ \dfrac{1}{\partial_k c_{\pm}(k_{m,0}^{\pm},\kappa)}\left(  \dfrac{1}{\alpha(k_{m,0}^{\pm},\kappa)} + \dfrac{3\ln 2}{\pi}  \right)\varepsilon + O(\varepsilon^2).
\end{equation}
\end{thm}

\medskip
Note that 
\begin{eqnarray*}
p&=&c_+(k,\kappa) \alpha(k,\kappa) + \dfrac{3\ln 2 \cdot \alpha}{\pi} \varepsilon + \varepsilon + O(s(\varepsilon)), \\
q&=&c_-(k,\kappa) \alpha(k,\kappa) + \dfrac{3\ln 2 \cdot \alpha}{\pi} \varepsilon + \varepsilon + O(t(\varepsilon)),
\end{eqnarray*}
and $k_{m,0}^{\pm}$ are roots of the leading-order terms $c_\pm(k,\kappa)=0$.
Hence the proof of the theorem follows the same perturbation argument as the one for Lemma 4.2 in the first part of this series \cite{lin_zhang17}, and we do not repeat it here. \\

\medskip

\noindent\textbf{Remark 3.1}  \textit{For a given $\kappa$, from Lemma  \ref{lem-disp_curve}, we have $|k_{m,0}^{\pm}(\kappa)|<|\kappa|$.
By assuming that $k_{m,0}$ is away from the cut-off frequency such that $k_{m,0}\in D_{\kappa,\delta, M}$, 
$|k_{m}^{\pm}(\kappa)|<|\kappa|$ holds true for $k_m^{\pm}$ obtained above. } \\

\noindent\textbf{Remark 3.2}  \textit{Since both $\partial_k c_{\pm}(k_{m,0}^{\pm},\kappa)$ and $\alpha(k_{m,0}^{\pm},\kappa)$ are real numbers, 
$k_{m,0}^{\pm}$ and the $O(\varepsilon)$ term in the above asymptotic expansion are real. In fact, since $|k_m^{\pm}|<|\kappa|$,
it can be argued  by variational method that $k_m^{\pm}$ are real eigenvalues. We refer to the Section 4.2 of \cite{lin_zhang17}
for a complete discussion. Therefore, the $O(\varepsilon^2)$ term in the asymptotic expansion is real too.}

\subsection{Asymptotic expansions of eigenmodes and surface bound states}\label{sec-H1_eigenmode}
For a given $\kappa$, recall that the eigenvectors for the corresponding two eigenvalues of $\mathbb{M}+\mathbb{I}$ 
are $[1 \quad 1]^T$ and $[1 \quad -1]^T$.
Therefore, if $k$ is an eigenvalues of the scattering operator such that $\lambda_1=0$ or $\lambda_2=0$, the solution of the
homogeneous linear system 
\begin{equation*}
( \mathbb{M}+\mathbb{I})\left[
 \begin{array}{llll}
\langle  \boldsymbol{\varphi}, \mathbf{e}_1 \rangle  \\
\langle  \boldsymbol{\varphi}, \mathbf{e}_2 \rangle
\end{array}
\right] =0
\end{equation*}
is given by
$$\left[
 \begin{array}{llll}
\langle  \boldsymbol{\varphi}, \mathbf{e}_1 \rangle  \\
\langle  \boldsymbol{\varphi}, \mathbf{e}_2 \rangle
\end{array}
\right] =c_1\left[
 \begin{array}{llll}
1 \\
1
\end{array}
\right] 
\quad\mbox{and}\quad
\left[
\begin{array}{llll}
\langle  \boldsymbol{\varphi}, \mathbf{e}_1 \rangle  \\
\langle  \boldsymbol{\varphi}, \mathbf{e}_2 \rangle
\end{array}
\right] =c_2\left[
 \begin{array}{llll}
1 \\
-1
\end{array}
\right] 
$$
respectively for some constant $c_1$ and $c_2$.

We derive the eigenmode of the homogeneous scattering problem. Without loss of generality, let us set $c_1=c_2=1$.
First consider the far-field zones $\Omega_1^+:=\{ x \;|\;  x_2 > 2 \}$ and $\Omega^-_{1}:=\{ x \;|\;  x_2 <-1 \}$ above and below the slab
respectively. By the quasi-periodicity of the solution, we may restrict the discussion 
to the domain $\Omega_1^+\cap\Omega^{(0)}$. Observe that the scattered field 
\begin{equation}\label{eq-u_far-field_exp}
u^s_\varepsilon(x) = \int_{\Gamma^+_\varepsilon} g^e(x,y) \dfrac{\partial u_\varepsilon(y)}{\partial \nu} ds_y = -\varepsilon \int_{0}^1 g^e(x,(\varepsilon Y, 1))  \varphi_1(Y)d Y \quad x\in\Omega_1^+\cap\Omega^{(0)}.
\end{equation}
Let $k=k_m^+$, then $\langle  \boldsymbol{\varphi}, \mathbf{e}_1 \rangle=1$. In addition,
\begin{equation}\label{eq-ge_far-field_asy}
g^e(x,(\varepsilon Y, 1)) = g^e(x,(0, 1))\left(1 + O(\varepsilon)\right) \quad \mbox{for} \; x\in \Omega^+ _{1}\cap\Omega^{(0)},
\end{equation}
we obtain
\begin{equation}\label{eq-u_far-field_asy}
u^s_\varepsilon(x) =  -\varepsilon  \left(1 + O(\varepsilon)\right)  g^e(x,(0, 1)) \quad \mbox{for} \; x\in \Omega^+ _{1}\cap\Omega^{(0)}.
\end{equation}

Note that $\kappa_n\sim O(1/\varepsilon)$  and
$\zeta_n(k)= \sqrt{k^2-\kappa_n^2}\sim O(1/\varepsilon)$ for $n\ne0$, since $d\sim\varepsilon$.
Therefore,
\begin{eqnarray}\label{eq-ge_far-field_exp}
g^e(x,(0, 1)) = 2g^d(x,(0,1)) &=& -\dfrac{i}{d}  \sum_{n=-\infty}^{\infty} \dfrac{1}{\zeta_n(k_m^+)   } e^{i \kappa_n(x_1)+i\zeta_n(k_m^+)   |x_2-1| } \nonumber \\
                          &=& -\dfrac{i}{d}  \dfrac{1}{\zeta_0(k_m^+)   } e^{i \kappa x_1+i\zeta_0(k)   |x_2-1| } +  O (e^{ -2\pi\eta/\varepsilon \cdot  |x_2-1| }) 
\end{eqnarray}
By substituting into \eqref{eq-u_far-field_asy} and using the fact that $|\kappa|>k_m^+$, it yields that
\begin{equation}\label{eq-mode1_H1}
u^s_\varepsilon(x) = \dfrac{\eta}{\sqrt{\kappa^2-(k_m^+)^2}   } e^{i \kappa x_1 -  \sqrt{\kappa^2-(k_m^+)^2}  \, |x_2-1| } + O(\varepsilon) \quad \mbox{for} \; x\in \Omega^+ _{1}\cap\Omega^{(0)}.
\end{equation}
Similarly, by using $\langle  \boldsymbol{\varphi}, \mathbf{e}_2 \rangle=1$, we have
\begin{equation}\label{eq-mode2_H1}
u^s_\varepsilon(x) = \dfrac{\eta}{\sqrt{\kappa^2-(k_m^+)^2}   } e^{i \kappa x_1 -  \sqrt{\kappa^2-(k_m^+)^2}  \, |x_2| } + O(\varepsilon) \quad \mbox{for} \; x\in \Omega^- _{1}\cap\Omega^{(0)}.
\end{equation}
Namely, the eigenmode is a surface bound-state mode that decays exponentially above and below the slab.
The same holds for eigenmode corresponding to $k=k_m^-$.

In the reference slit $S_{\varepsilon}^{(0)}$, by noting that
\begin{equation*} 
\left\{
\begin{array}{llll}
\vspace*{0.1cm}
\Delta u_{\varepsilon} + k^2 u_{\varepsilon} = 0 \quad \mbox{in} \; S_{\varepsilon}^{(0)},  \\
\vspace*{0.1cm}
\dfrac{\partial u_{\varepsilon}}{\partial x_1} = 0  \quad \mbox{on} \; x_1=0, \,\, x_1=\varepsilon,
\end{array}
\right.
\end{equation*}
we may expand $u_\varepsilon$ as the sum of wave-guide modes as follows:
\begin{equation}\label{eq-u_slit}
u_\varepsilon(x)= a_0 e^{ikx_2} + b_0 e^{ik(1-x_2)} + \sum_{m\geq 1} \left(a_m  e^{-k_2^{(m)}x_2} + b_m e^{-k_2^{(m)}(1-x_2)} \right) \cos \dfrac{m\pi x_1}{\varepsilon} ,
\end{equation}
where $k_2^{(m)}=\sqrt{(m\pi/\varepsilon)^2-k^2}$.
Taking the derivative of \eqref{eq-u_slit} and evaluating on the slit apertures, it follows that
\begin{eqnarray}\label{eq-du_expansion}
\dfrac{\partial u_\varepsilon}{\partial x_2}(x_1,1)&=& ika_0 e^{ik}  - ikb_0 + \sum_{m\geq 1} \left( -a_m e^{-k_2^{(m)}} + b_m 
\right) k_2^{(m)} \cos \dfrac{m\pi x_1}{\varepsilon}, \label{eq-du_expansion1}  \\
\dfrac{\partial u_\varepsilon}{\partial x_2}(x_1,0)&=&
 ika_0 - ikb_0 e^{ik}+\sum_{m\geq 1} \left( -a_m  + b_m e^{-k_2^{(m)}}
\right) k_2^{(m)} \cos \dfrac{m\pi x_1}{\varepsilon} \label{eq-du_expansion2}.
\end{eqnarray}
For $k=k_m^+$, recall that
$$ \langle  \boldsymbol{\varphi}, \mathbf{e}_1 \rangle= \langle  \boldsymbol{\varphi}, \mathbf{e}_2 \rangle=1, $$
and from \eqref{Op_eqns_1},
$$   \boldsymbol{\varphi} =  -(\beta+\tilde\beta) (\mathbb{L}^{-1}  \mathbf{e}_1 +  \mathbb{L}^{-1}  \mathbf{e}_2). $$
Therefore, it can be shown that 
$$ a_0  = a_0^+ = \dfrac{e^{-ik}}{(1- e^{ik})}, \quad b_0  = b_0^+ = \dfrac{e^{-ik}}{(1- e^{ik})}, $$
and 
$$ \abs{a_m} \le C/\sqrt{m} , \quad \abs{b_m} \le C/\sqrt{m}, \quad \mbox{for}\; \,m\ge1. $$
A similar calculation for $k=k_m^+$ leads to
$$ a_0 = a_0^- = -\dfrac{e^{-ik}}{(1+e^{ik})}, \quad b_0  = b_0^- = -\dfrac{e^{ik}}{(1+e^{ik})}. $$
Therefore, for a given $\kappa$, the eigenmode in the slit region $S_{\varepsilon}^{(0), int} :=\{ x\in S_{\varepsilon}^{(0)} \;|\; x_2 \gg \varepsilon , 1- x_2 \gg \varepsilon \} $ adopts the following asymptotical expansion:
$$
u_\varepsilon(x) =a_0^\pm e^{ikx_2} + b_0^\pm e^{ik(1-x_2)}  + O\left(e^{-1/\varepsilon} \right) .
$$
for the eigenvalue $k=k_m^\pm$.

\subsection{Homogenization and effective medium theory}\label{sec-H1_effective_medium}
As $\varepsilon\to0$, by the homogenization theory, one expects that the scattering by the slab with an array of slits 
is equivalent to the scattering by a homogenous slab medium. To this end, let us consider the incident wave 
$u^{i}= e^{i ( \kappa x_1 - \zeta (x_2-1))}$ that impinges on the slab, where $\kappa=k\sin\theta$ and $\zeta=k\cos\theta$.
The calculations for $u_\varepsilon$ in the far-field zone are parallel to the ones presented in Section \ref{sec-H1_eigenmode}.  
First, it is clear that the scattered field $u_\varepsilon^s$ is given by \eqref{eq-u_far-field_exp} in $\Omega^+ _{1}\cap\Omega^{(0)}$ . 
Using the asymptotic expansion of the Green's function \eqref{eq-ge_far-field_asy}, it follows that
\begin{equation*}
u^s_\varepsilon(x) =  -\varepsilon  \left(1 + O(\varepsilon)\right) \cdot  g^e(x,(0, 1))  \cdot  \int_{0}^1  \varphi_1(Y)d Y  \quad \mbox{for} \; x\in \Omega^+ _{1}\cap\Omega^{(0)}.
\end{equation*}
An application of the asymptotic expansions for the Green's function in \eqref{eq-ge_far-field_exp} and $\langle  \boldsymbol{\varphi}, \mathbf{e}_1 \rangle$
in Lemma \ref{lem-phi} leads to
\begin{eqnarray}\label{eq-far_field_scat}
u^s_\varepsilon(x) &=&  -\varepsilon  \left(1 + O(\varepsilon)\right) \cdot  \left( -\dfrac{i}{d}  \dfrac{1}{\zeta_0(k)   } e^{i \kappa x_1+i\zeta_0(k)   |x_2-1| } +  O (e^{ -1/\varepsilon}) \right)  \cdot  \bigg(\alpha+ O(r(\varepsilon))  \bigg) \left(\dfrac{1}{p} + \dfrac{1}{q}\right) \nonumber \\
&=& \dfrac{i\varepsilon \alpha  }{d \; \zeta} \cdot \left(\dfrac{1}{p} + \dfrac{1}{q}\right) \cdot  e^{i (\kappa x_1+\zeta (x_2-1)) }  \cdot  \left(1 + O(\varepsilon)\right).
\end{eqnarray}
Therefore, by virtue of Lemma \ref{eq-u_eps_formula} and the relation $\varepsilon=\eta d$,  the total field
$$ u_\varepsilon(x) = u^i(x) + \left[1+ \dfrac{i\eta \alpha  }{\zeta} \cdot \left(\dfrac{1}{p} + \dfrac{1}{q}\right) \cdot  \left(1 + O(\varepsilon)\right) \right] \cdot e^{i (\kappa x_1+\zeta (x_2-1))}   \quad \mbox{for} \; x\in \Omega^+ _{1}\cap\Omega^{(0)}. $$
A straightforward calculation based on explicit expressions \eqref{eq-p_H1} and \eqref{eq-q_H1} gives
$$  \alpha \cdot \left(\dfrac{1}{p} + \dfrac{1}{q}\right) = \dfrac{2\zeta^2k-i\cdot2\zeta\eta \cdot k^2\tan k}{ -(\zeta^2 + \eta^2 k^2)\tan k - i \cdot 2 \zeta \eta k  } \left(1 + O(\varepsilon)\right) . $$
We substitute the above into \eqref{eq-far_field_scat} and obtain
\begin{equation}\label{eq-far_field_up}
u_\varepsilon(x) = u^i(x) + R \cdot e^{i (\kappa x_1+\zeta (x_2-1))}   \quad \mbox{for} \; x\in \Omega^+ _{1}\cap\Omega^{(0)},
\end{equation}
where the reflection coefficient
\begin{equation}\label{eq-R}
R = \dfrac{i\cdot(-\zeta^2+\eta^2k^2)\tan k}{ -i\cdot (\zeta^2 + \eta^2 k^2)\tan k + 2 \zeta \eta k  } \cdot \left(1 + O(\varepsilon)\right).  
\end{equation}
Similarly, it can be obtained that the transmitted field below the slab is
\begin{equation}\label{eq-far_field_bottom}
u_\varepsilon(x) = T \cdot e^{i (\kappa x_1-\zeta x_2)}   \quad \mbox{for} \; x\in \Omega^- _{1}\cap\Omega^{(0)},
\end{equation}
where the transmission coefficient
\begin{equation}\label{eq-T}
T = \dfrac{2\zeta \eta k}{ -i\cdot (\zeta^2 + \eta^2 k^2)\sin k + 2 \zeta \eta k \cos k  } \cdot \left(1 + O(\varepsilon)\right). 
\end{equation}
\medskip

\begin{figure}[!htbp]
\begin{center}
\includegraphics[height=4.5cm,width=12cm]{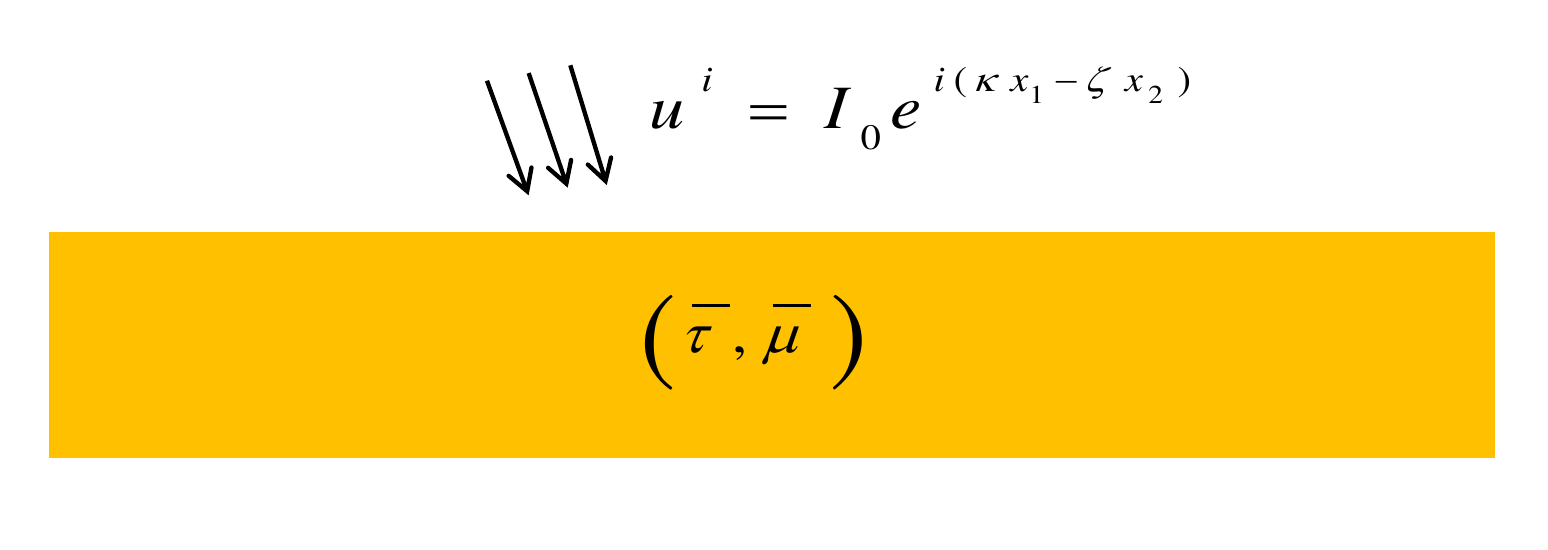}
\vspace*{-20pt}
\caption{The effective layered medium as $\varepsilon\to0$.}
\label{fig-eff_med}
\vspace{-10pt}
\end{center}
\end{figure}

Now let us derive the effective slab medium as $\varepsilon\to0$. 
Denote the relative permittivity and the permeability of the effective medium in the slab by $\bar{\tau}$ and $\bar{\mu}$ respectively,
and consider the layered medium as depicted in Figure \ref{fig-eff_med}. The corresponding scattering problem is formulated
as 
\begin{equation}\label{eq-scattering_hom}
\nabla\cdot \left(\dfrac{1}{\tau} \nabla u\right) + k^2 \mu u =0,
\end{equation}
where
\begin{equation*}
\tau(x_1,x_2)  = \left\{
\begin{array}{lll}
\vspace*{5pt}
1,  & x_2>1 \; \mbox{or} \; x_2<0, \\
\bar{\tau},  & 0<x_2<1. \\
\end{array}
\right.
\quad \mbox{and} \quad
\mu(x_1,x_2)  = \left\{
\begin{array}{lll}
\vspace*{5pt}
1,  & x_2>1 \; \mbox{or} \; x_2<0, \\
\bar{\mu},  & 0<x_2<1. \\
\end{array}
\right.
\end{equation*}

We look for $\bar{\tau}$ and $\bar{\mu}$ such that the associated far-field $u$ recovers
the leading-order term of the far-field $u_\varepsilon$ given by \eqref{eq-far_field_up} - \eqref{eq-T}.
\begin{thm} \label{thm-homo1}
Let
$$
\bar\tau=\left[
 \begin{array}{llll}
 \infty & 0 \\
 0   & 1/\eta
\end{array}
\right] \quad \mbox{and}  \quad\bar\mu=\eta. $$
If the incident wave $u^{i}= I_0 e^{i ( \kappa x_1 - \zeta (x_2-1))}$, where $\kappa=k\sin\theta$ and $\zeta=k\cos\theta$,
then the total field for the scattering problem \eqref{eq-scattering_hom} has the following form
\begin{equation*}
u(x_1,x_2)  = \left\{
\begin{array}{lll}
\vspace*{5pt}
I_0 e^{i (\kappa x_1 - \zeta (x_2-1)} +   R_0 e^{i (\kappa x_1 + \zeta (x_2-1))} ,  & x_2>1,  \\
T_0 e^{i (\kappa x_1 - \zeta x_2)},  & x_2<0
\end{array}
\right.
\end{equation*}
The reflection and transmission coefficients are given by
\begin{equation*}
R_0 =  \dfrac{i\cdot(-\zeta^2+\eta^2k^2)\tan k}{ -i\cdot (\zeta^2 + \eta^2 k^2)\tan k + 2 \zeta \eta k }\cdot I_0 \quad \mbox{and} \quad
T_0 = \dfrac{2\zeta \eta k}{ -i\cdot (\zeta^2 + \eta^2 k^2)\sin k + 2 \zeta \eta k \cos k  }\cdot I_0.
\end{equation*}
\end{thm}

\noindent\textbf{Proof}.
If $\bar\tau$ and $\bar\mu$ are given as in the theorem, then the solution of the scattering problem can be written down as
follows in each layer:
\begin{equation*}
u(x_1,x_2)  = \left\{
\begin{array}{lll}
\vspace*{5pt}
 I_0 e^{i (\kappa x_1 - \zeta (x_2-1))} +   R_0 e^{i (\kappa x_1 + \zeta (x_2-1))} ,  & x_2>1,  \\
 a^+  e^{i (\kappa x_1 + k x_2)} +   a^-e^{i (\kappa x_1 - k x_2)} ,  & 0<x_2<1, \\
 T_0 e^{i (\kappa x_1 - \zeta x_2)},  & x_2<0.
\end{array}
\right.
\end{equation*}
By imposing the continuity conditions along the interfaces $x_2=0$ and $x_2=1$:
\begin{eqnarray*}
&& u(x_1,0-) = u(x_1, 0+), \quad \partial_{x_2} u(x_1,0-) = \eta\partial_{x_2}u(x_1, 0+); \\
&& u(x_1,1-) = u(x_1, 1+), \quad \partial_{x_2} u(x_1,1-) = \eta\partial_{x_2}u(x_1, 1+),
\end{eqnarray*}
we obtain the following linear system for $(R_0, T_0, a^+, a^-)$:
\begin{eqnarray*}
a^+ + a^- &=& T_0, \hspace*{3.6cm}   -i\zeta T_0 = i \eta k (a^+ - a^-),\\
e^{ik}a^+ + c^{-ik} a^- &=& I_0+R_0, \quad i \eta k (e^{ik}a^+ - e^{-ik} a^-) = i\zeta (-I_0+ R_0 ).
\end{eqnarray*}
This can be further reduced to the following system:
\begin{eqnarray}\label{eq-linear_system_h1}
\left[-i \tan k \cdot (1 + \tilde\zeta^2) + 2 \tilde\zeta \right]  R_0 &=& I_0 \cdot i \cdot (1-\tilde\zeta^2) \tan k, \label{eq-linear_system_r0} \\
\left[-i \tan k \cdot (1 + \tilde\zeta^2) + 2 \tilde\zeta \right]  T_0 &=& I_0 \cdot 2 \tilde\zeta /\cos k,\label{eq-linear_system_t0} \\
\frac{1}{2}\left( 1-\tilde\zeta\right) T_0 - a^+ &=& 0 ,  \label{eq-linear_system_c+} \\
\frac{1}{2}\left( 1+\tilde\zeta \right) T_0 - a^- &=& 0,  \label{eq-linear_system_c-}
\end{eqnarray}
where $\tilde\zeta=\zeta/(\eta k)$. Solving \eqref{eq-linear_system_r0} and \eqref{eq-linear_system_t0} proves the assertion. \qed \\

Next, we demonstrate the dispersion relation for the homogenized layered medium recovers the leading-order term 
of the dispersion relation $k_m^\pm(\kappa)$ given in Theorem \ref{thm-asym_dispersion}. 
\begin{thm} \label{thm-homo2}
If 
$$
\bar\tau=\left[
 \begin{array}{llll}
 \infty & 0 \\
 0   & 1/\eta
\end{array}
\right] \quad \mbox{and}  \quad\bar\mu=\eta, $$
then the dispersion relation for the layered medium have two branches given by
\begin{equation}\label{eq-dispersion_hom_med}
\kappa=k \sqrt{1+ \eta^2 \tan^2(k/2)} \quad \mbox{and} \quad \kappa=k \sqrt{1+ \eta^2 \cot^2(k/2)}.
\end{equation}
The corresponding eigenmode is
\begin{equation*}
u(x_1,x_2)  = \left\{
\begin{array}{lll}
\vspace*{5pt}
R_0 e^{i \kappa x_1 +  \sqrt{\kappa^2-k^2} x_2} ,  & x_2>1,  \\
 a^+  e^{i (\kappa x_1 + k x_2)} +   a^-e^{i (\kappa x_1 - k x_2)} ,   & 0<x_2<1, \\
 e^{i \kappa x_1 - \sqrt{\kappa^2-k^2} x_2},  & x_2<0.
\end{array}
\right.
\end{equation*}
where $R_0= \dfrac{i\sin k}{2}  \cdot \left(\dfrac{\eta k}{\zeta} - \dfrac{\zeta}{\eta k}\right), a^+= \dfrac{1}{2}\left( 1-\dfrac{\zeta}{\eta k} \right), a^- = \dfrac{1}{2}\left( 1+\dfrac{\zeta}{\eta k} \right).
$
\end{thm}

\noindent\textbf{Proof}. 
To obtain the dispersion relation, we solve for $(k, \kappa)$ such that there exists nontrivial solutions for the linear system \eqref{eq-linear_system_r0} - \eqref{eq-linear_system_c-} when $I_0=0$.
This implies that 
$$ 
-i \tan k \cdot (1 + \tilde\zeta^2) + 2 \tilde\zeta=0
$$
so the determinant of the coefficient matrix is zero. 
Solving the above equation yields
$$ \tilde\zeta = \dfrac{\zeta}{\eta k} = i \tan(k/2) \quad \mbox{or} \quad \tilde\zeta=\dfrac{\zeta}{\eta k} = -i \cot(k/2). $$
Using the relation $\kappa^2+\zeta^2 = k^2$, it follows that
$$ \kappa^2=k^2 (1+ \eta^2 \tan^2(k/2)) \quad \mbox{or} \quad  \kappa^2=k^2 (1+ \eta^2 \cot^2(k/2)).$$
Finally, the corresponding nontrivial solutions to the above linear system are
$$
T_0=C, \quad R_0= \frac{i e^{-i\zeta}}{2} \sin k \cdot \left(\frac{\eta k}{\zeta} - \frac{\zeta}{\eta k}\right) C,  \quad
 a^+= \dfrac{1}{2}\left( 1-\dfrac{\zeta}{\eta k} \right) C,  \quad a^- = \dfrac{1}{2}\left( 1+\dfrac{\zeta}{\eta k} \right) C
$$
for some constant $C$. By taking $C=1$, we 
proved the second part of the theorem. \qed
\\

\subsection{Surface plasmon for plasmonic metals and perfect conductors with slits}\label{sec-H1_spoof_spp}
\begin{figure}[!htbp]
\begin{center}
\includegraphics[height=3.5cm,width=9cm]{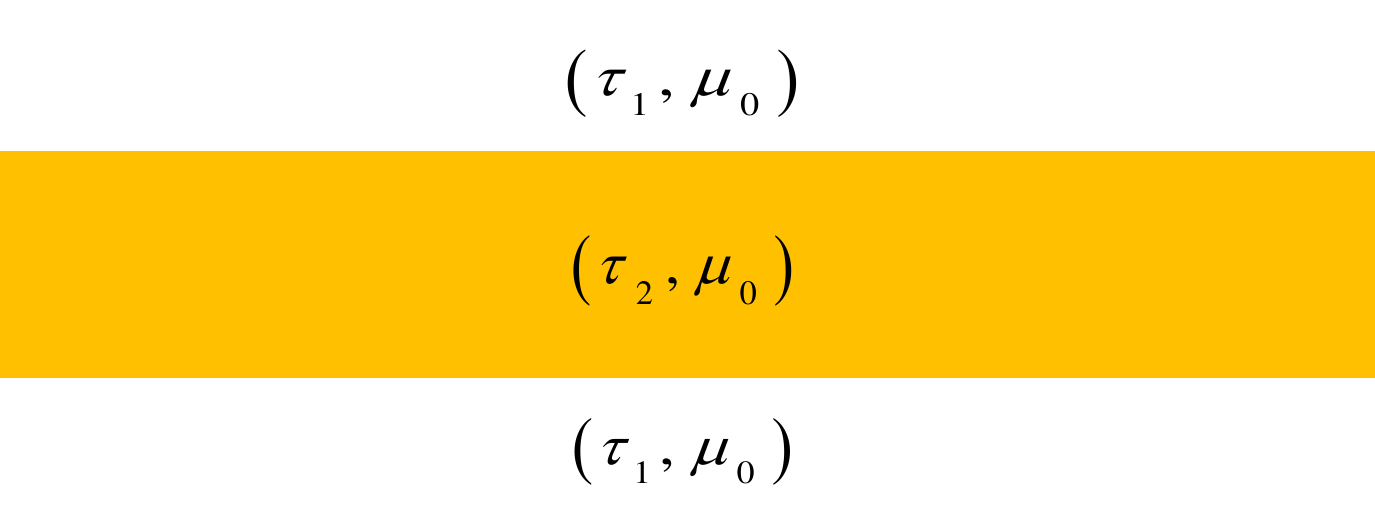}\hspace*{-0.5cm}
\includegraphics[height=4.5cm,width=7.5cm]{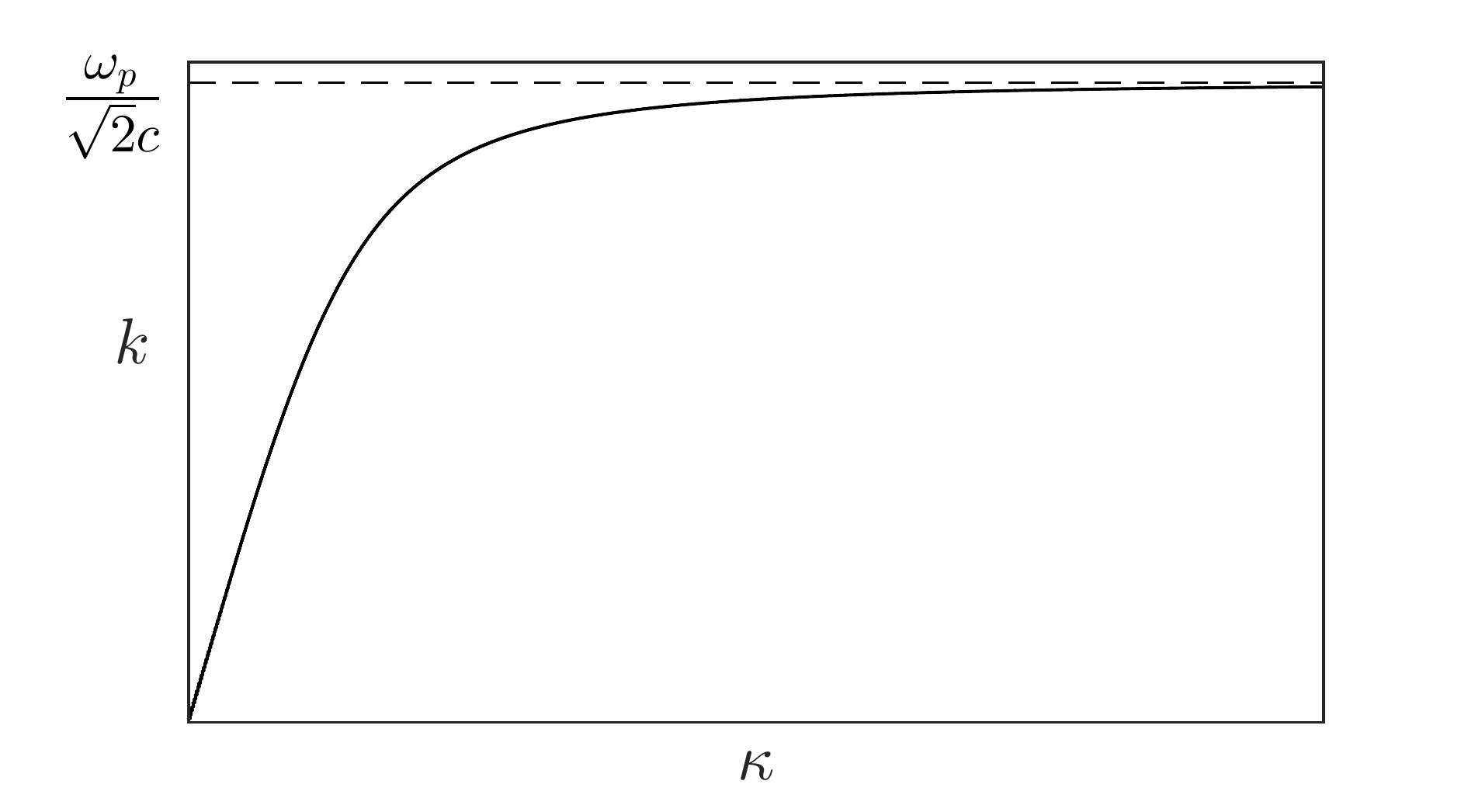}
\vspace{-10pt}
\caption{The diaelctric-metal-diaelectric medium (left) and the associated dispersion curve (right).}\label{fig-dispersion_metal_slab}
\vspace{-10pt}
\end{center}
\end{figure}

It is known that surface plasmon modes are supported on the flat interface of dielectric and noble metal. 
Let the permittivity of the dielectric material and the metal be $\tau_1$ and $\tau_2$ respectively, and $\Re \; \tau_2<0$.
Then it can be calculated that, for a metal slab with a thickness of $\ell$,
the following localized modes exist along the interfaces of the dielectric-metal medium  (\cite{maier07})
\begin{equation}\label{eq-mode-metal_slab}
u(x_1,x_2)  = \left\{
\begin{array}{lll}
\vspace*{5pt}
 e^{i \kappa x_1 -  \sqrt{\kappa^2-k^2\tau_1} \; x_2} ,  & x_2>0,  \\
 a^+  e^{i (\kappa x_1 + \sqrt{\kappa^2-k^2\tau_2} \; x_2)} +   a^-e^{i (\kappa x_1 - \sqrt{\kappa^2-k^2\tau_2} \; x_2)} ,   & 0<x_2<\ell, \\
 t_0 e^{i \kappa x_1 + \sqrt{\kappa^2-k^2} x_2},  & x_2<0.
\end{array}
\right.
\end{equation}
In addition, the dispersion relations are given by
$$ \tanh ( \sqrt{\kappa^2-k^2\tau_2} \ell)  + \dfrac{\tau_1\sqrt{\kappa^2-k^2\tau_2}}{\tau_2\sqrt{\kappa^2-k^2\tau_1}} =0 
\quad \mbox{and}  \quad
\tanh ( \sqrt{\kappa^2-k^2\tau_2} \ell)  + \dfrac{\tau_2\sqrt{\kappa^2-k^2\tau_1}}{\tau_1\sqrt{\kappa^2-k^2\tau_2}} =0.
$$
For simplicity, assume that the exterior medium is vacuum so that $\tau_1=1$.
If one applies the Drude model without damping for the metal permittivity by letting $\tau_2=1-\dfrac{\omega_p^2}{\omega^2}$, where
$\omega_p$ is the plasma frequency and it takes the value $\omega_p=1.37\times10^{16}$Hz for gold \cite{ordal83}, 
then the first dispersion relation is shown in Figure \ref{fig-dispersion_metal_slab}, and
the second dispersion relation has a similar shape.

A direct comparison of Figure \ref{H1_dispersion} and \ref{fig-dispersion_metal_slab} confirms the resemblance of the dispersion curves
for the plasmonic metal and the perfect conductor with slits as $\varepsilon\to0$. Both dispersion curves lie below the light line
such that $k(\kappa)<|\kappa|$, and approach certain freququencies as $\kappa\to\infty$.
In addition, from \eqref{eq-mode1_H1},  \eqref{eq-mode2_H1} and \eqref{eq-mode-metal_slab}, the corresponding eigenmodes are both localized 
bound states along the slab interfaces.
That is, surface plasmonic effect mimicking that of plasmonic metals exists in a perfect conducting slab by engineering its surfaces.
In particular, for a PEC slab with a thickness of $\ell$, by a scaling argument,
it is seen that the wavenumber $k \to \pi/\ell$ and $2\pi/\ell$ respectively for the first branch of two dispersion curves as the $\kappa$ increases to infinity.
While for the plasmonic metal, the wavenumber $k \to \omega_p/(\sqrt{2} c)$ as $\kappa$ increases to infinity. 
Therefore, $1/\ell$ determines the plasmonic frequency for the perfect conductor. As such one can tune the associated plasmonic mode in different frequencies by adjusting the thickness of the metal slab $\ell$.

\subsection{Total transmission for the scattering by an incident plane wave}\label{sec-H1_total_trans}
As discussed in previous sections, surface bound states occur when $k^\pm(\kappa)<|\kappa|$.
Now if one considers scattering by an incident plane wave $u^{i}= e^{i ( \kappa x_1 - \zeta (x_2-1))}$, where $\kappa=k\sin\theta$ and $\zeta=k\cos\theta$.
Then $|\kappa|<k$ holds, and the solution to the scattering problem is unique. The corresponding reflection and transmission coefficients are
given by \eqref{eq-R} and \eqref{eq-T}. As $\varepsilon\to0$, their limit values are the ones associated with the effective medium as stated in 
Theorem \ref{thm-homo1}. In this section, we investigate the field pattern above and below the metal slab in the limiting case of $\varepsilon\to0$.
To this end, let us rewrite the reflection coefficient $R_0$  and the transmission coefficient $T_0$ in Theorem \ref{thm-homo1} as
\begin{eqnarray*}
R_0 &=&  \dfrac{i  \tan k \cdot (\eta^2- \cos^2 \theta)}{-i \tan k \cdot (\eta^2+ \cos^2 \theta)) + 2 \eta \cos \theta}, \\
T_0 &=& \dfrac{2 \cos \theta \cdot \eta }{-i \sin k \cdot (\eta^2+ \cos^2 \theta) + 2 \cos \theta \cdot \eta \cos k }.
\end{eqnarray*}

\begin{figure}[!htbp]
\begin{center}
\includegraphics[height=4.5cm,width=8cm]{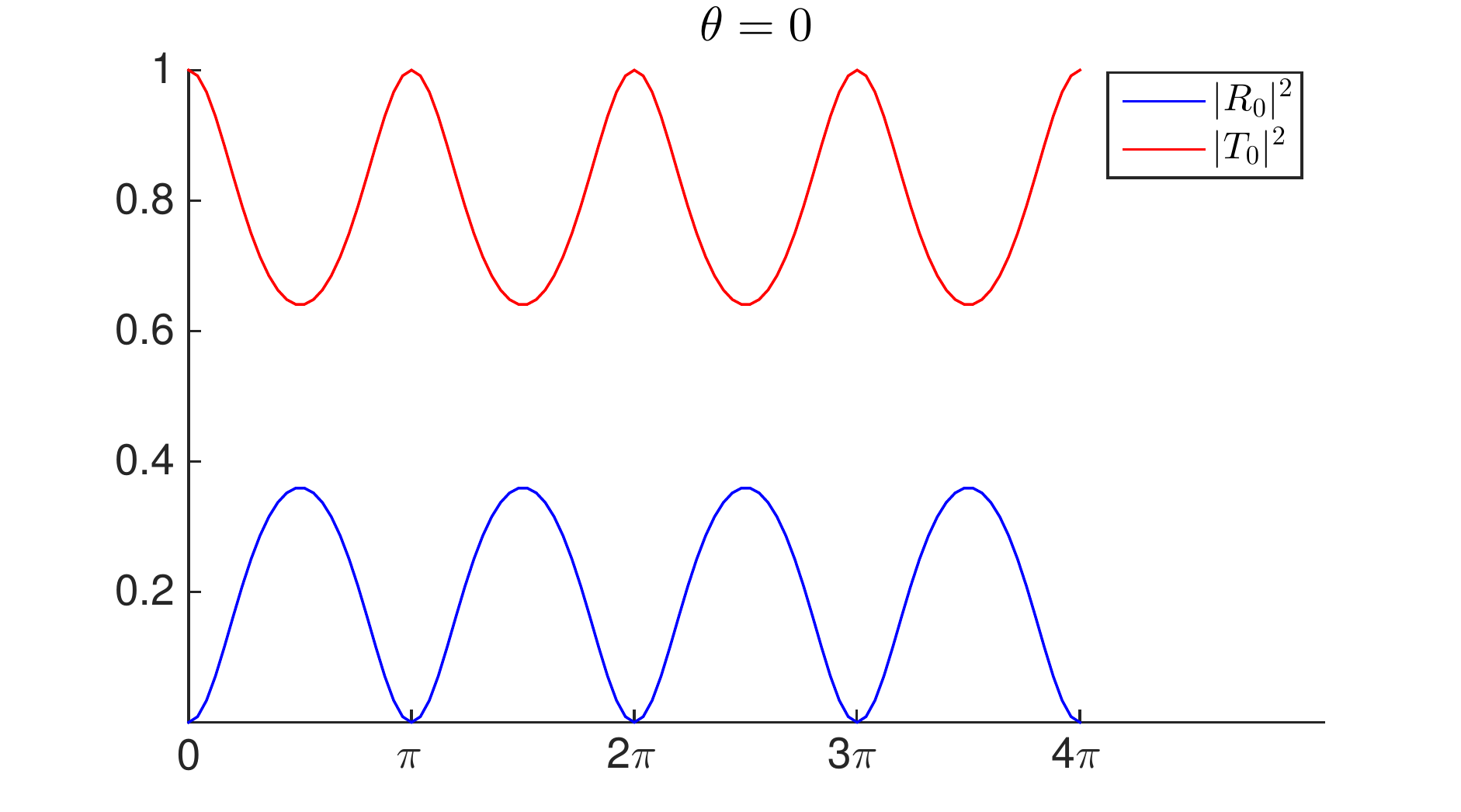}
\includegraphics[height=4.5cm,width=8cm]{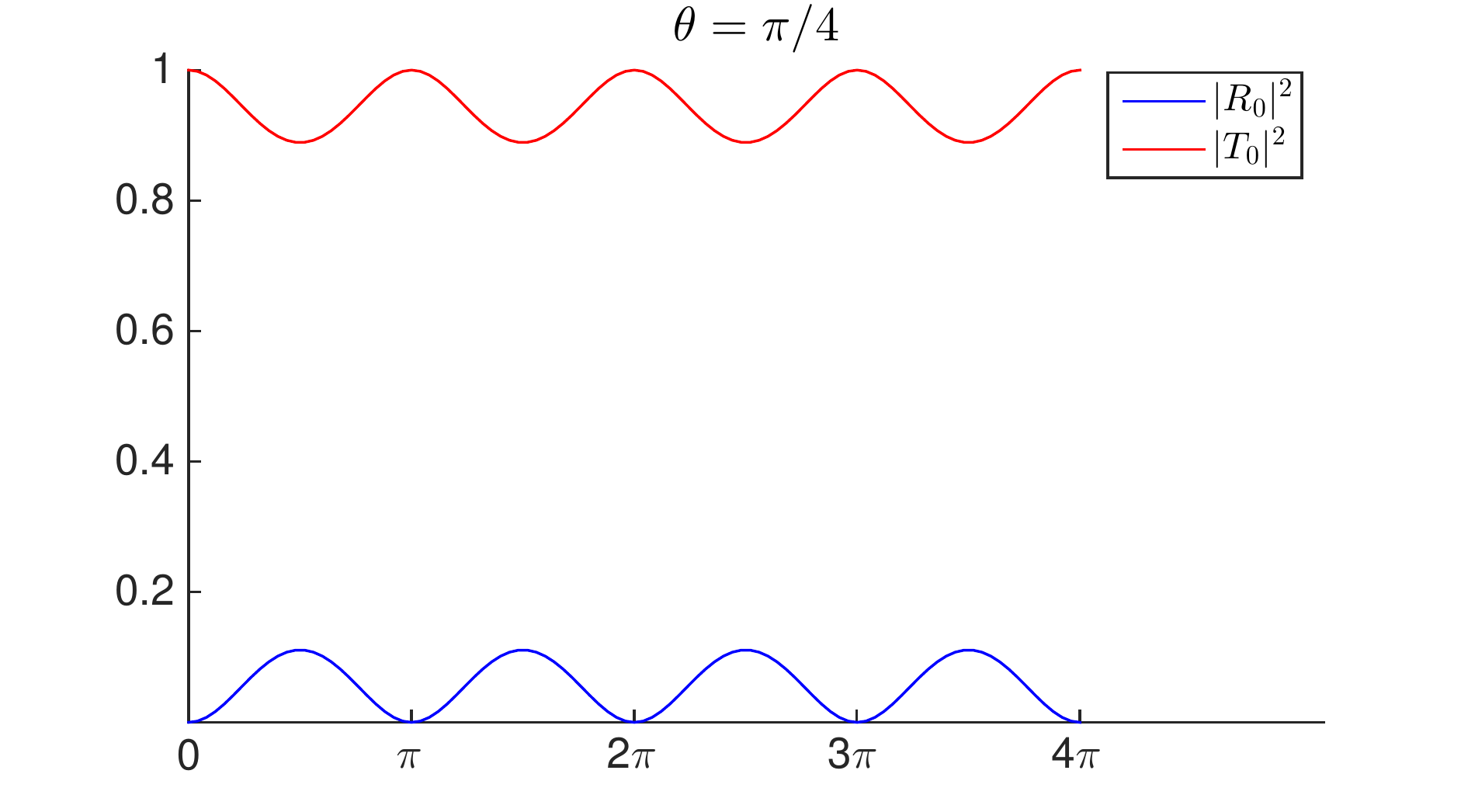}
\includegraphics[height=4.5cm,width=8cm]{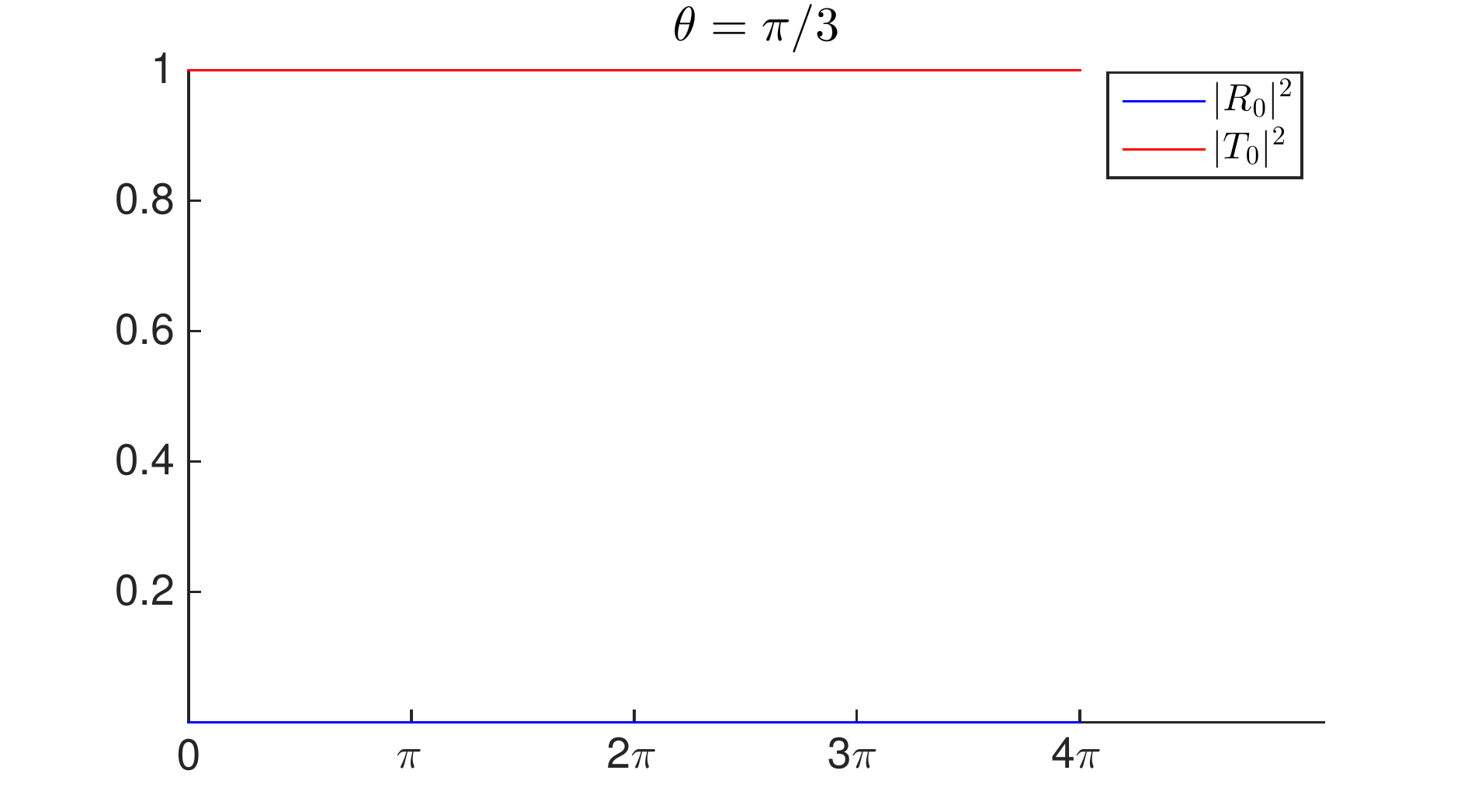}
\includegraphics[height=4.5cm,width=8cm]{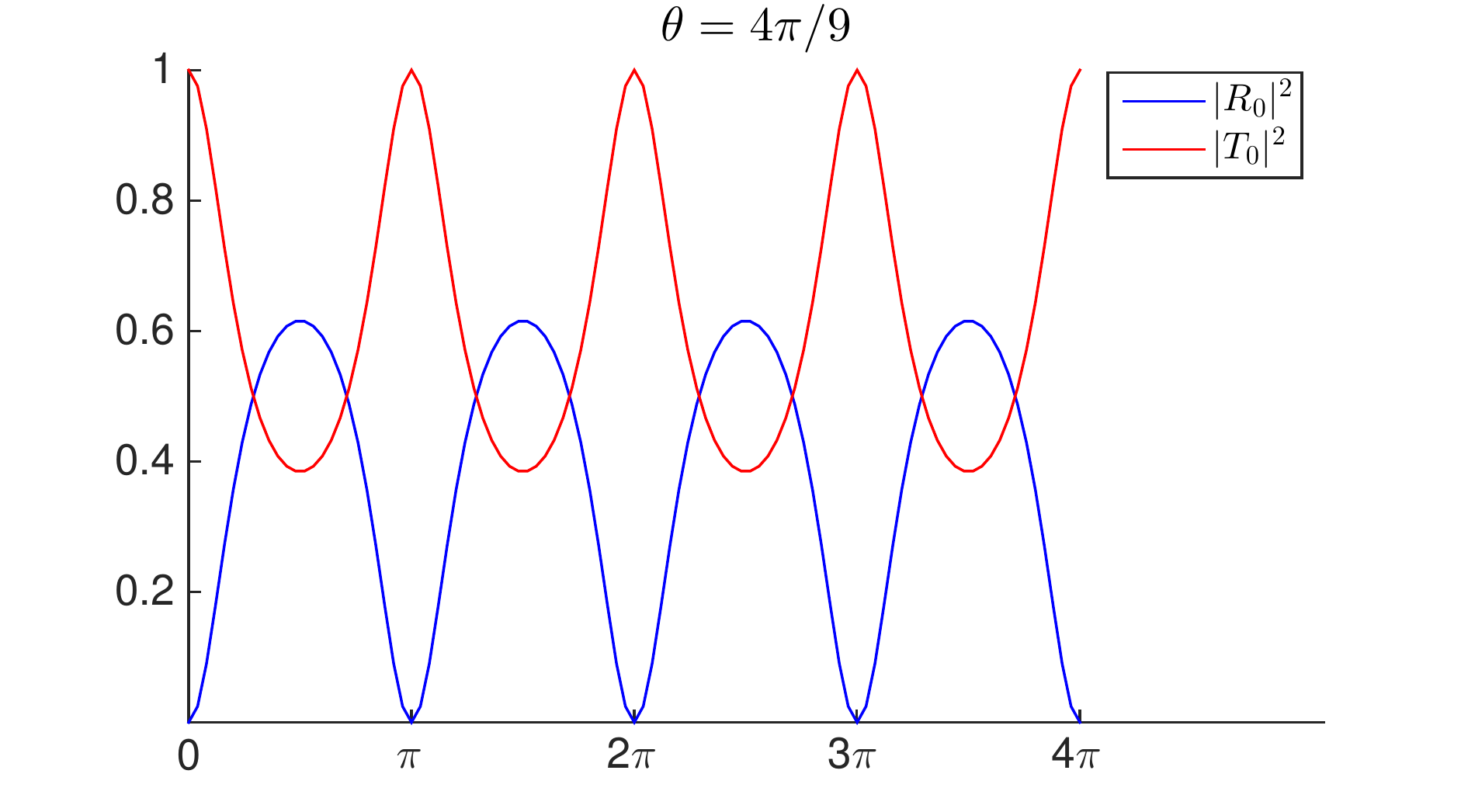}
\vspace{-10pt}
\caption{$|R_0|^2$ and $|T_0|^2$ for various incident angles and wavenumbers when $\eta=0.5$. Note that $|R_0|^2+|T_0|^2=1$.}\label{fig-total_trans}
\vspace{-10pt}
\end{center}
\end{figure}

When $\eta=0.5$, their amplitudes for various incident angles and wavenumbers are shown in Figure \ref{fig-total_trans}.
It is seen that when $k=m\pi$, where $m$ is an integer, $|T_0|=1$ for all incident angles. 
That is, total transmission is achieved at those wavenumbers by the scattering of the homogenized slab, which is viewed as the limiting
effective medium of a perfect conducting slab perorated with an array of small slits and with small periods.
For the special incident angle such that $\cos\theta=\eta$, total transmission is obtained throughout all the frequencies (see Figure \ref{fig-total_trans}).
We note that perfect transmission has also been reported for highly conductive metals patterned with narrow slits \cite{bouchitte-schweizer13}.

Since $k^\pm(\kappa)<|\kappa|$ holds for the real dispersion curves \eqref{eq-dispersion_hom_med},
thus for the incident plane wave with $\kappa=k\sin\theta$ and $|\kappa|<k$, the frequencies
$k=m\pi$ are not associated with ``plasmonic frequencies" given by the dispersion relation,
or certain scattering resonances which are defined as the poles of the resolvent associated with the scattering problem.
Furthermore,  based on Lemma \ref{lem-u_slit} which is given in Section \ref{sec-H2}, and the asymptotic expansion of $p$ and $q$,
it can be shown that no field enhancement occurs inside the slits at those wavenumbers when $\varepsilon$ is small
(see Section \ref{sec-H2_EM_period} for a discussion when $k\to0$). 
Hence, we deduce that the total transmission observed here is not due to plamonic resonant effect or scattering resonance.
Instead, it may be due to the so-called Fabry-Perot resonances associated with the homogenized slab in Section \ref{sec-H1_effective_medium}, for which all reflected waves from the slab boundaries interfere destructively 
and zero reflected wave is finally attained on top of the slab \cite{vaughan89}. On the other hand, for the total transmission at the special incident angle such that $\cos\theta=\eta$, the physical mechanism is not quite clear.

\section{Homogenization regime (H2): non-resonant field enhancement} \label{sec-H2}
In the homogenization regime (H2) where $\varepsilon \ll d \ll \lambda$ (see Figure \ref{fig-geo_H1H2}, bottom), 
there exists no resonance or eigenvalue such that the homogeneous scattering problem 
attains nontrivial solutions. Namely, the corresponding scattering problem \eqref{eq-Helmholtz}-\eqref{eq-rad_cond} attains a unique solution. This is demonstrated in Section \ref{sec-H2_res_eig}. 
In Section \ref{sec-H2_EM_fields}, we derive the asymptotic expansion of the wave fields in both the near and far field zones, and study
their enhancement behaviors in this regime. It is shown that although no enhancement is gained for the magnetic field,
strong electric field is induced in the slits and on the slit apertures.
A discussion on the field enhancement for varying period $d$ is presented in Section \ref{sec-H2_EM_period}.
Briefly speaking, the field enhancement becomes stronger as $d$ increases. For extremely large $d$ that still satisfies $d \ll \lambda$, the effect of periodicity is vanishing and enhancement behavior resembles that of the single slit considered in \cite{lin_zhang16} as $d\to\infty$.
On the other hand, as $d$ decreases, the field enhancement becomes weaker. In particular, if $d\sim\varepsilon$ holds,
then no electromagnetic field enhancement is gained.

\subsection{Non-existence of resonance or eigenvalue}\label{sec-H2_res_eig}
From \eqref{eq-scattering4}, the homogeneous scattering problem with the incident wave $u^i=0$
can be equivalently formulated as the operator equation
$$ (\mathbb{P} + \mathbb{L})\boldsymbol{\varphi}=0, $$
which further reduces to
\begin{equation*}
( \mathbb{M}+\mathbb{I})\left[
 \begin{array}{llll}
\langle  \boldsymbol{\varphi}, \mathbf{e}_1 \rangle  \\
\langle  \boldsymbol{\varphi}, \mathbf{e}_2 \rangle
\end{array}
\right] =0,
\end{equation*}
by \eqref{eq-linear_sys}.
Therefore, the resonances/eigenvalues of the scattering operator are roots of $\lambda_1(k; \kappa, d, \varepsilon) $ and $\lambda_2(k; \kappa, d, \varepsilon)$, the eigenvalues of $\mathbb{M}+\mathbb{I}$.
Equivalently, they are roots of $p(k; \kappa, d, \varepsilon) =0$ and $q(k; \kappa, d, \varepsilon) =0$.

Let us define
\begin{equation}\label{eq-gamma}
\gamma(k, \kappa, d)=\dfrac{1}{\pi} \left(3\ln 2 + \ln\dfrac{\pi}{d}\right) + \left(\dfrac{1}{2\pi}\sum_{n\neq 0} \dfrac{1}{|n|} - \dfrac{i}{d} \sum_{n=-\infty}^{\infty}  \dfrac{1}{\zeta_n(k) } \right),
\end{equation}
where 
\begin{equation*}
\zeta_n(k) =\zeta(k, \kappa, d)  = \left\{
\begin{array}{lll}
\vspace*{5pt}
\sqrt{k^2-(\kappa+2\pi n/d)^2},  & \abs{\kappa+2\pi n/d} < k, \\
i\sqrt{(\kappa+2\pi n/d)^2-k^2},  & \abs{\kappa+2\pi n/d} > k. \\
\end{array}
\right.
\end{equation*}
Then from \eqref{eq-lambda1}, the definition of $\beta$ in \eqref{beta}, and Lemma \ref{lem-L_inv}, we may explicitly express
\begin{eqnarray}\label{formula_p_H2}
p(k; \kappa, d, \varepsilon) &=& \varepsilon +\left[ \dfrac{\cot k }{k} + \dfrac{1}{k\sin k} +  \varepsilon \gamma(k, \kappa, d) + \dfrac{1}{\pi} \varepsilon \ln \varepsilon   \right]   \left(\langle \mathbb{L}^{-1}  \mathbf{e}_1,  \mathbf{e}_1 \rangle  + \langle \mathbb{L}^{-1}  \mathbf{e}_1,  \mathbf{e}_2\rangle\right)  \nonumber \\
&=& \varepsilon +\left[ \dfrac{\cot k }{k} + \dfrac{1}{k\sin k} +  \varepsilon \gamma(k, \kappa, d) + \dfrac{1}{\pi} \varepsilon \ln \varepsilon   \right]   \left(\alpha + s(\varepsilon) \right),
\end{eqnarray}
where $s(\varepsilon)\sim O(r(\varepsilon))$. Similarly,
\begin{equation}\label{formula_q_H2}
q(k; \kappa, d, \varepsilon)=\varepsilon +\left[ \dfrac{\cot k }{k} - \dfrac{1}{k\sin k} +  \varepsilon \gamma(k, \kappa, d) + \dfrac{1}{\pi} \varepsilon \ln \varepsilon   \right]   \left(\alpha + t(\varepsilon) \right),
\end{equation}
where $t(\varepsilon)\sim O(r(\varepsilon))$. It is clear that as $k\to0$, the leading order of $p$ and $q$
$$  \left(\dfrac{\cot k }{k} + \dfrac{1}{k\sin k}\right) \alpha \to \infty \quad \mbox{and} \quad  
\left(\dfrac{\cot k }{k} - \dfrac{1}{k\sin k}\right) \alpha \to -\dfrac{\alpha}{2} $$
respectively. Therefore, as $\varepsilon\to0$, $p(k; \kappa, d, \varepsilon) =0$ and $q(k; \kappa, d, \varepsilon) =0$ do not attain roots for $k \ll 1$ .

\subsection{Quantitative analysis of the electromagnetic field in the near-field and far-field zones}\label{sec-H2_EM_fields}
\subsubsection{Field enhancement in the slits}
From previous discussion, the scattering problem in the homogenization regime (H2) attains a unique solution.
In this section, we investigate the electromagnetic field in both near-field and far-field zones.
Note that in the reference slit $S_{\varepsilon}^{(0)}$,  from Section \ref{sec-H1_eigenmode},
$u_\varepsilon$ can be expanded as
\begin{equation}\label{u_in_slit}
u_\varepsilon(x)= a_0 \cos kx_2 + b_0 \cos k(1-x_2) + \sum_{m\geq 1} \left(a_m e^{-k_2^{(m)}x_2} +e^{-k_2^{(m)}x_2}  \right) \cos \dfrac{m\pi x_1}{\varepsilon}, 
\end{equation}
where $k_2^{(m)}=\sqrt{(m\pi/\varepsilon)^2-k^2}$.
The following asymptotic expansion holds for $u_\varepsilon$ in  $S_{\varepsilon}^{(0), int} :=\{ x\in S_{\varepsilon}^{(0)} \;|\; x_2 \gg \varepsilon , 1- x_2 \gg \varepsilon \} $:

 \begin{lem}\label{lem-u_slit}
In the slit region $S_{\varepsilon}^{(0), int}$, we have $u_\varepsilon(x_1,x_2)=u_0(x_2)+ u_\infty(x_1,x_2)$, where
\begin{equation}\label{eq-u_0}
u_0(x_2) = \bigg[ \alpha+ O(r(\varepsilon)) \bigg] \left[ \dfrac{\cos (kx_2)}{k\sin k}  \left(\dfrac{1}{p} + \dfrac{1}{q}\right)
+ \dfrac{\cos (k(1-x_2)) }{k\sin k }   \left(\dfrac{1}{p} - \dfrac{1}{q}\right) \right],
\end{equation}
and $u_\infty\sim O\left(e^{-1/\varepsilon}\right)$. Here $\alpha$ is defined in Lemma \ref{lem-optK}.
 \end{lem}
 \noindent\textbf{Proof}
 From the expansion \eqref{u_in_slit}, it follows that
\begin{eqnarray}\label{eq-du_expansion}
\dfrac{\partial u_\varepsilon}{\partial x_2}(x_1,1)&=& ika_0 e^{ik}  - ikb_0 + \sum_{m\geq 1} \left( -a_m e^{-k_2^{(m)}} + b_m 
\right) k_2^{(m)} \cos \dfrac{m\pi x_1}{\varepsilon}, \label{eq-du_expansion1}  \\
\dfrac{\partial u_\varepsilon}{\partial x_2}(x_1,0)&=&
 ika_0 - ikb_0 e^{ik}+\sum_{m\geq 1} \left( -a_m  + b_m e^{-k_2^{(m)}}
\right) k_2^{(m)} \cos \dfrac{m\pi x_1}{\varepsilon} \label{eq-du_expansion2}.
\end{eqnarray}
Therefore,
\begin{eqnarray*}
- a_0 k\sin k  &=& \dfrac{1}{\varepsilon} \int_{\Gamma^+_\varepsilon} \dfrac{\partial u_\varepsilon}{\partial x_2}(x_1,1) dx_1 
= -\int_{0}^1 \varphi_1(X)d X = -\bigg[\alpha+ O(r(\varepsilon))  \bigg] \left(\dfrac{1}{p} + \dfrac{1}{q}\right), \\
b_0 k \sin k &=& \dfrac{1}{\varepsilon} \int_{\Gamma^-_\varepsilon} \dfrac{\partial u_\varepsilon}{\partial x_2}(x_1,0) dx_1 
= \int_{0}^1 \varphi_2(X)d X = \bigg[\alpha+ O(r(\varepsilon)) \bigg] \left(\dfrac{1}{p} - \dfrac{1}{q}\right).
 \end{eqnarray*}
We obtain
\begin{equation}\label{a0b0}
a_0  = \dfrac{1}{k \sin k}\bigg[\alpha+ O(r(\varepsilon))  \bigg] \left(\dfrac{1}{p} + \dfrac{1}{q}\right),  \quad
b_0 = \dfrac{1}{k \sin k} \bigg[ \alpha+ O(r(\varepsilon)) \bigg]  \left(\dfrac{1}{p} - \dfrac{1}{q}\right).  
\end{equation}

For $m\ge1$, the coefficients $a_m$ and $b_m$ can be obtained similarly by taking the inner product of \eqref{eq-du_expansion1}
and \eqref{eq-du_expansion2} with $\cos \dfrac{m\pi x_1}{\varepsilon}$. Then a direct estimate leads to
\begin{equation}\label{ambm}
\abs{a_m} \le C/\sqrt{m} , \quad \abs{b_m} \le C/\sqrt{m}, \quad \mbox{for}\; \,m\ge1,
\end{equation}
where $C$ is some positive constant independent of $\varepsilon$, $k$ and $m$.
The proof is complete by substituting \eqref{a0b0} and \eqref{ambm} into \eqref{u_in_slit}. \qed \\

Recall that in homogenization regime (H2), $\varepsilon \ll 1$ and $k\ll 1$ holds. In what follows, we set $k=\varepsilon^\sigma$, where $\sigma>0$.
\begin{lem}\label{lem-pq_H2}
Let $\sigma>0$ and $k=\varepsilon^\sigma$, then 
\begin{equation}\label{eq-p_H2}
\dfrac{1}{p} \cdot \dfrac{1}{k\sin k} = \dfrac{1}{ 2 \alpha}  (1+O(\varepsilon^{2\sigma})+O(\varepsilon^{\sigma+1})),
\end{equation}
and
\begin{equation}\label{eq-q_H2}
\dfrac{1}{q} = \left\{
\begin{array}{lll}
-\dfrac{2}{\alpha}\left(1+O(\varepsilon^{2\sigma}) + O(\varepsilon^{1-\sigma})\right)  & \,\, \mbox{if} \,\,\, 0<\sigma<1,\\ \\
\dfrac{i\cdot d \cos\theta}{\alpha}\varepsilon^{\sigma-1} \left(1+O(\varepsilon^{\sigma-1})\right)  & \,\, \mbox{if} \,\, \, \sigma>1,
\end{array}
\right.
\end{equation}
\end{lem}
\noindent\textbf{Proof} From the expression of $\gamma$ in \eqref{eq-gamma}, it is clear that
$$ \gamma(k, \kappa, d) = -\dfrac{i}{d\;\zeta_0(k)} + O(1) = -\dfrac{i}{k d \cos\theta} + O(1) $$
if $k\ll 1$.
From the explicit formulas of $p$ and $q$ in \eqref{formula_p_H2} and  \eqref{formula_q_H2}, a direct calculation yields
\begin{eqnarray}\label{eq-p_H2_der}
\dfrac{1}{p} \cdot \dfrac{1}{k\sin k} &=& \dfrac{1}{ (\cos k +1) \alpha} \left(1 + O(k^2\varepsilon \ln\varepsilon) + O(\gamma k^2\varepsilon) \right) \nonumber\\
&=&\dfrac{1}{ 2 \alpha (1+O(\varepsilon^{2\sigma}))} \left(1 + O(\varepsilon^{2+2\sigma} \ln\varepsilon) + \frac{1}{d} \cdot O(\varepsilon^{1+\sigma}) \right) \nonumber \\
&=&\dfrac{1}{ 2 \alpha}  \left(1+O(\varepsilon^{2\sigma}) +\frac{1}{d} \cdot  O(\varepsilon^{1+\sigma})  \right).
\end{eqnarray}
On the other hand,
\begin{eqnarray}\label{eq-q_H2_der}
q &=& \left( -\dfrac{1}{2}+O(k^2) + \gamma\varepsilon+ \frac{1}{\pi}\varepsilon \ln\varepsilon \right)\big(\alpha+t(\varepsilon)\big)+\varepsilon  \nonumber \\
   &=&  -\alpha\left(\dfrac{1}{2}+O(\varepsilon^{2\sigma}) + \frac{i}{d \cos\theta}\varepsilon^{1-\sigma}+ O(\varepsilon \ln\varepsilon)\right),
\end{eqnarray}
whence the asymptotic expansion of $1/q$ follows. \qed \\

\begin{thm}\label{thm-u_slit}
Let $\sigma>0$ and $k=\varepsilon^\sigma$, then $u_\varepsilon(x_1,x_2)=u_0(x_2)+ u_\infty(x_1,x_2)$, where
$$
 u_0(x_2) = \left\{
\begin{array}{lll}
2x_2 +O(\varepsilon^{2\sigma}) + O(\varepsilon^{1-\sigma})  &  \mbox{if}\,\, \; 0<\sigma<1,\\ \\
1+ i d \cdot\cos\theta  \big(2x_2-1\big) \varepsilon^{\sigma-1}  +O(\varepsilon^{\sigma+1})  + O(\varepsilon^{2(\sigma-1)}) &  \mbox{if} \,\, \;  \sigma>1,
\end{array}
\right.
$$
and  $u_\infty\sim O\left(e^{-1/\varepsilon}\right)$.
\end{thm}

 \noindent\textbf{Proof} By a combination of Lemma \ref{lem-u_slit} and \ref{lem-pq_H2}, and the Taylor expansion, it follows that when $0<\sigma<1$,
\begin{eqnarray*}
u_0(x_2) &=& \bigg( 1+ O(r(\varepsilon)) \bigg)\bigg[ \frac{1}{2} \big(1+O(\varepsilon^{2\sigma})\big) \big(\cos (kx_2) + \cos(k(1-x_2))\big)  \\  
&& -2\big(1+O(\varepsilon^{2\sigma}) + O(\varepsilon^{1-\sigma})\big)  \dfrac{\cos (kx_2) - \cos(k(1-x_2))}{k\sin k} \bigg] \\
&=& \bigg( 1+ O(r(\varepsilon)) \bigg)\bigg[1+O(\varepsilon^{2\sigma})-\left(1+O(\varepsilon^{2\sigma}\big) + O(\varepsilon^{1-\sigma})\right)   
\big(1-2x_2\big) \bigg]  \\
&=& 2x_2 +O(\varepsilon^{2\sigma}) + O(\varepsilon^{1-\sigma}).
\end{eqnarray*}
While for $\sigma>1$,
\begin{eqnarray*}
u_0(x_2) &=& \bigg( 1+ O(r(\varepsilon)) \bigg)\bigg[1+O(\varepsilon^{\sigma+1})+\frac{id \cdot \cos\theta}{2} \; \varepsilon^{\sigma-1} \left(1+O(\varepsilon^{\sigma-1})\right)\big(1-2x_2\big) \bigg]  \\
&=& 1+ \frac{i d \cdot\cos\theta}{2}  \big(1-2x_2\big) \varepsilon^{\sigma-1} +O(\varepsilon^{\sigma+1}) + O(\varepsilon^{2(\sigma-1)}).
\end{eqnarray*}

\qed

From the above theorem, we see that there is no enhancement for the magnetic field $u_\varepsilon$ in the homogenization regime (H2).
However, the transition of the magnetic field $u_\varepsilon$ along the $x_2$ direction resembles 
a linear function with a slope of $2$ (for $0< \sigma <1$) and  $O(\varepsilon^{\sigma-1})$ (for $\sigma >1$) in the slits. This is in contrast with the incident field, which changes with a rate of 
$O(k)$, or $O(\varepsilon^{\sigma})$, in the slits.
Such fast transition of magnetic field from the upper to lower slit aperture, compared to the incident wave, induces strong electric field enhancement as stated in the following theorem.
\begin{thm}\label{thm-E_slit}
If $\varepsilon \ll 1$ and $k=\varepsilon^\sigma$, 
then the electric field $E_{\varepsilon}=[E_{\varepsilon,1}, E_{\varepsilon,2}, 0]$ in $S_{\varepsilon}^{int}$, where
$$ 
E_{\varepsilon,1} = \left\{
\begin{array}{lll} 
\dfrac{2i}{\varepsilon^\sigma \sqrt{\tau_0/\mu_0}}+ \min\{O(\varepsilon^{\sigma}), O(\varepsilon^{1-2\sigma}) \}&  \mbox{if} \; 0<\sigma<1,\\ \\
\dfrac{d \cos\theta}{\varepsilon \sqrt{\tau_0/\mu_0}} + \min\{O(\varepsilon), O(\varepsilon^{\sigma-2}) \} &  \mbox{if} \;  \sigma>1,
\end{array}
\right.
\quad \mbox{and} \quad
E_{\varepsilon,2}  \sim O(e^{-1/\varepsilon}/\varepsilon^\sigma ). $$
$\tau_0$ and $\mu_0$ is the electric permittivity and magnetic permeability in the vacuum respectively.
\end{thm}
\noindent\textbf{Proof} Note that in the TM case, the magnetic field is given by
$$
H_\varepsilon=[0, 0, u_\varepsilon]. 
$$
Therefore, by Ampere's law 
$$\nabla\times  H_\varepsilon= [ \partial u_\varepsilon/\partial x_2,  -\partial u_\varepsilon/\partial x_1, 0]= -i\omega\tau_0 E_\varepsilon.$$
For $0<\sigma<1$, we have
\begin{align*} 
E_{\varepsilon,1} &= \dfrac{2i}{k \sqrt{\tau_0/\mu_0}} +  O(\varepsilon^{2\sigma}/k) +  O(\varepsilon^{1-\sigma}/k) = \dfrac{2i}{\varepsilon^\sigma \sqrt{\tau_0/\mu_0}} 
+  O(\varepsilon^{\sigma}) +  O(\varepsilon^{1-2\sigma}),\\
E_{\varepsilon,2} &= -\partial u_\infty/\partial x_1\cdot i/ \omega\tau_0    \sim O(e^{-1/\varepsilon}/k) = O(e^{-1/\varepsilon}/\varepsilon^\sigma ). 
\end{align*}
The electric field when $\sigma>1$ follows by a similar calcuation.
\qed

\begin{rmk}
From the above theorem, we see that the enhancement for the electric field is not uniform throughout the low frequency regime.
 When $k=\varepsilon^\sigma$ and $0<\sigma<1$,  $E_\varepsilon$ is of order $O(1/\varepsilon^\sigma)$,
or equivalently $O(1/k)$. Thus the enhancement becomes stronger as $k$ decreases in such scenario. While for $\sigma>1$,  $E_\varepsilon$ is of order $O(1/\varepsilon)$, which is independent of $k$.
\end{rmk}

\begin{rmk}
It is also observed from the previous discussion that the electric field enhancement also depends on the size of period $d$.
Such dependence is significant when $\sigma>1$. This will be discussed in more details in Section \ref{sec-H2_EM_period}.
\end{rmk}

\subsubsection{Field enhancement on apertures of slits}
Define
\begin{equation} \label{eq-h}  
h(X) = \dfrac{1}{\pi}\int_{0}^1 \ln | X-Y| (K^{-1}1)(Y)dY,
\end{equation}
and let
\begin{equation}\label{beta_bar}
\bar \beta^e(k,\kappa, d) := \beta^e(k,\kappa, d, \varepsilon)-\dfrac{1}{\pi} \ln \varepsilon  =\left(\ln 2 + \ln\dfrac{\pi}{d}\right) + \left(\dfrac{1}{2\pi}\sum_{n\neq 0} \dfrac{1}{|n|} - \dfrac{i}{d} \sum_{n=-\infty}^{\infty}  \dfrac{1}{\zeta_n(k)}\right).
\end{equation}

\begin{lem}\label{lem-u_aperture}
The following asymptotic holds for the total field
\begin{eqnarray}\label{u_up_aperture}
u_\varepsilon(x_1, 1)  &=& -\dfrac{1}{\pi}  \left( 
\dfrac{\alpha}{p} +\dfrac{\alpha}{q}  \right) \cdot  \varepsilon \ln\varepsilon - 
\left( \dfrac{\alpha}{p} + \dfrac{\alpha}{q}  \right) \left(\bar\beta^e +h(x_1/\varepsilon) \right) \cdot  \varepsilon  + 2  \nonumber \\
&& - \left(\dfrac{\alpha}{p} + \dfrac{\alpha}{q}\right) \cdot O(\varepsilon\ln\varepsilon \cdot r(\varepsilon))-\kappa \cdot O(\varepsilon)+O(\varepsilon \cdot r(\varepsilon)) 
\end{eqnarray}
and
\begin{eqnarray}\label{u_low_aperture}
u_\varepsilon(x_1, 0)  &=& -\dfrac{1}{\pi}  \left( 
\dfrac{\alpha}{p} -\dfrac{\alpha}{q}  \right) \cdot  \varepsilon \ln\varepsilon - 
\left( \dfrac{\alpha}{p} - \dfrac{\alpha}{q}  \right) \left(\bar\beta^e +h(x_1/\varepsilon) \right) \cdot  \varepsilon  \nonumber \\
&& - \left(\dfrac{\alpha}{p} - \dfrac{\alpha}{q}\right) \cdot O(\varepsilon\ln\varepsilon \cdot r(\varepsilon))+O(\varepsilon \cdot r(\varepsilon))
\end{eqnarray}
on the slit apertures $\Gamma^+_\varepsilon$ and $\Gamma^-_\varepsilon$ respectively. 
\end{lem}

 \noindent\textbf{Proof} Recall that on $\Gamma^+_\varepsilon$,
$$ u_\varepsilon(x) = \int_{\Gamma^+_\varepsilon} g_\varepsilon^e(x,y) \dfrac{\partial u_\varepsilon (y)}{\partial \nu} ds_y + u^i+ u^r.$$
Let $x_1= \varepsilon X $, $y_1= \varepsilon Y$.
We have 
\begin{equation*}
u_\varepsilon(\varepsilon X, 1) =  -\int_{0}^1 G_\varepsilon^e(X,Y)\varepsilon \varphi_1(Y)dY +f(X). 
\end{equation*} 
Using Lemma \ref{lem-phi} and the asymptotic expansion of $G_\varepsilon^e(X,Y)$ in Lemma \ref{lem-periodic_green}, we obtain
\begin{eqnarray*}
u_\varepsilon(\varepsilon X, 1) &=& -\varepsilon \beta^e\bigg(\alpha+ O(r(\varepsilon)) \bigg) \left(\dfrac{1}{p} + \dfrac{1}{q}\right) 
 - \dfrac{\varepsilon}{\pi} \left( \kappa \cdot O(1) + \dfrac{\alpha}{p} + \dfrac{\alpha}{q} \right)\int_{0}^1 \ln | X-Y| (K^{-1}1)(Y)dY \\
&& - \left(\dfrac{\alpha}{p} + \dfrac{\alpha}{q} \right) O(\varepsilon \cdot r(\varepsilon))  +  O(\varepsilon \cdot r(\varepsilon)) + f(X).
 \end{eqnarray*}
The desired expansion follows by using \eqref{eq-h} and \eqref{beta_bar}. The wave field on the lower aperture can be obtained similarly. \qed \\
 
Now if $\sigma>0$ and $k=\varepsilon^\sigma$, By subsituting \eqref{eq-p_H2}-\eqref{eq-q_H2}
into the above lemma,  it follows that
\begin{equation*}
u(x_1,1)=2+O(\varepsilon\ln\varepsilon),  \quad u(x_1,0)=O(\varepsilon\ln\varepsilon),
\end{equation*}
and there is no enhancement for the magnetic field on the aperture. 
The enhancement of the electric field is stated in the following Theorem.
\begin{thm}\label{thm-E_aperture}
Let $\sigma >0$ and $k=\varepsilon^\sigma$,  then the following hold for electric field
$$ 
E_\varepsilon(x_1,1)= \left\{
\begin{array}{lll} 
 \dfrac{2i}{\varepsilon^\sigma\sqrt{\tau_0/\mu_0}} [  K^{-1}1,  -h'(X), 0] +\min\{O(\varepsilon^{\sigma}), O(\varepsilon^{1-2\sigma}) \} &  \mbox{if}\,\, \; 0<\sigma<1,\\ \\
\dfrac{d \cos\theta}{\varepsilon \sqrt{\tau_0/\mu_0}} [  K^{-1}1,  -h'(X), 0] +   O(\varepsilon^{\sigma-2}) &  \mbox{if}\,\, \;  \sigma>1,
\end{array}
\right.
$$ 
$$ 
E_\varepsilon(x_1,0)= \left\{
\begin{array}{lll} 
 \dfrac{2i}{\varepsilon^\sigma\sqrt{\tau_0/\mu_0}} [  K^{-1}1,  h'(X), 0] +\min\{O(\varepsilon^{\sigma}), O(\varepsilon^{1-2\sigma}) \} &  \mbox{if}\,\, \; 0<\sigma<1,\\ \\
\dfrac{d \cos\theta}{\varepsilon \sqrt{\tau_0/\mu_0}} [  K^{-1}1,  h'(X), 0] + O(\varepsilon^{\sigma-2})  &  \mbox{if} \,\,\;  \sigma>1,
\end{array}
\right.
$$ 
on the upper and lower apertures respectively.
\end{thm}
\noindent\textbf{Proof}
We derive $E_\varepsilon$ on the upper slit apertures. The case for the lower slit apertures can be obtained similarly.
Taking the derivative of \eqref{u_up_aperture} yields
\begin{eqnarray*}
\dfrac{\partial u_\varepsilon}{\partial x_1}(x_1,1)  &=&   -\left( \dfrac{\alpha}{p} +\dfrac{\alpha}{q} \right) \cdot \dfrac{1}{\varepsilon} h'(X)  \cdot  \varepsilon  - \left(\dfrac{\alpha}{p} + \dfrac{\alpha}{q}\right) \cdot O(\varepsilon\ln\varepsilon \cdot r(\varepsilon))-\kappa \cdot O(\varepsilon),
\end{eqnarray*}
where  $h(X)$ is defined by (\ref{eq-h}).
Therefore, using \eqref{eq-p_H2}-\eqref{eq-q_H2}, we see that
\begin{equation}\label{du_dx1_up_H2}
\frac{\partial u_\varepsilon}{\partial x_1}(x_1,1)  =
\left\{
\begin{array}{lll} 
 2  h'(X) + O(\varepsilon^{2\sigma}) + O(\varepsilon^{1-\sigma}),  &  \mbox{if} \,\, \; 0<\sigma<1, \\ \\
 -i\cdot d \cos\theta \; h'(X) \cdot  \varepsilon^{\sigma-1} + O(\varepsilon^{2(\sigma-1)})  &  \mbox{if}\,\, \; \sigma>1.
\end{array}
\right.
\end{equation}
On the other hand, by \eqref{eq-varphi} it follows that
\begin{eqnarray*}
\dfrac{\partial u_\varepsilon}{\partial x_2}(x_1,1) &=&  -K^{-1}1  \cdot \left( \kappa \cdot O(\varepsilon) + \dfrac{\alpha}{p}  + \dfrac{\alpha}{q}  \right) + \left(\dfrac{\alpha}{p} +  \dfrac{\alpha}{q} \right) \cdot O(r(\varepsilon)) +  O(r(\varepsilon)).
\end{eqnarray*}
An application of \eqref{eq-p_H2}-\eqref{eq-q_H2} yields
\begin{equation}\label{du_dx2_aperture_H2}
\frac{\partial u_\varepsilon}{\partial x_2}(x_1,1) =
\left\{
\begin{array}{lll} 
 2 K^{-1}1 + O(\varepsilon^{2\sigma}) + O(\varepsilon^{1-\sigma}),  &  \mbox{if}\,\, \; 0<\sigma<1, \\ \\
 -i \cdot d \cos\theta \; K^{-1}1 \cdot \varepsilon^{\sigma-1} + O(\varepsilon^{2(\sigma-1)})  &  \mbox{if} \,\,\; \sigma>1.
\end{array}
\right.
\end{equation}
A combination of \eqref{du_dx1_up_H2} - \eqref{du_dx2_aperture_H2} and the Ampere's law leads to the desired asymptotic expansions
for the electric field on the apertures. \qed

\subsubsection{Far field asymptotic and effective medium theory}
In the far-field zone $\Omega_1^+:=\{ x \;|\;  x_2 > 2 \}$ above the slits, by restricting the domain to the reference period $\Omega_1^+\cap\Omega^{(0)}$,
we note that the scattered field
\begin{equation*}
u^s_\varepsilon(x) = \int_{\Gamma^+_\varepsilon} g^e(x,y) \dfrac{\partial u_\varepsilon(y)}{\partial \nu} ds_y.
\end{equation*}
An application of formula \eqref{du_dx2_aperture_H2} yields
$$u^s_\varepsilon \sim O(\varepsilon) \quad \mbox{and} \quad u^s_\varepsilon \sim O(\varepsilon^{2\sigma-1}) $$ 
for $0<\sigma<1$ and $\sigma>1$ respectively. The same holds true for the  far-field zone below the slits.
This shows that there is no electric or magnetic field enhancement in the far field. 
Moreover, as $\varepsilon\to0$, the effect of the slits vanishes and the perforated perfect conducting slab becomes a homogeneous perfect conducting slab.

\subsection{Electric field enhancement in the near field for varying sizes of period}\label{sec-H2_EM_period}
From Theorem \ref{thm-E_slit} and \ref{thm-E_aperture}, it is observed that the enhancement for the electric field  $E_\varepsilon$ depends on the size of the period $d$.
More precisely, if $k=\varepsilon^\sigma$, then for $0<\sigma<1$, the enhancement is of order 
$O(1/\varepsilon^\sigma)$ (or equivalently $O(1/k)$) and  is slightly affected as $d$ increases,
since $d$ appears in the high-order terms of $E_\varepsilon$. While for $\sigma>1$,  $d$ appears in the leading-order term of $E_\varepsilon$. In particular, the enhancement becomes stronger as $d$ increases.
Let us set $d=O(\varepsilon^{1-\sigma-\delta})$ for some $0<\delta<1$, then $d\ll\lambda$ still holds and $d\to\infty$ as $\varepsilon\to0$ in such scenario.
By substituting $d$ into \eqref{eq-p_H2_der} and  \eqref{eq-q_H2_der}, it is clear that the following lemma holds for $p$ and $q$.
\begin{lem}\label{lem-pq_large_H2}
If $\varepsilon \ll 1$, $k=\varepsilon^\sigma$ with $\sigma>1$, and  $d=O(\varepsilon^{1-\sigma-\delta})$ with $0<\delta<1$, then 
\begin{equation*}
\dfrac{1}{p} \cdot \dfrac{1}{k\sin k} = \dfrac{1}{ 2 \alpha}  (1+O(\varepsilon^{2\sigma})),
\end{equation*}
and
\begin{equation*}
\dfrac{1}{q} =-\dfrac{2}{\alpha}\left(1+O(\varepsilon^{2\sigma}) + O(\varepsilon^{\delta})+O(\varepsilon\ln\varepsilon)\right)
\end{equation*}
\end{lem}
Following the same lines as in Theorem \ref{thm-u_slit} and  \ref{thm-E_slit}, it can be shown that, for $\sigma>1$,
$$u= 2x_2 + O(\varepsilon^{\delta})
\quad\mbox{and}\quad 
E_{\varepsilon,1} =  \dfrac{2i}{\varepsilon^\sigma \sqrt{\tau_0/\mu_0}} + O(\varepsilon^{\delta-\sigma})$$
in the slits. Therefore, we recover the $O(1/\varepsilon^\sigma)$ order (or equivalently $O(1/k)$ order )  enhancement for $\sigma>1$.
Namely, for sufficiently large $d$, an uniform $O(1/k)$ enhancement for $E_\varepsilon$ is achieved throughout the low frequency regime.
This is consistent with the field enhancement for a single slit perorated in a perfect conducting slab (when $d=\infty$),
where an enhancement order of $O(1/k)$ is obtained  throughout the low frequency regime \cite{lin_zhang16}.

One the other hand,  
as the period $d$ decreases, the magnitude of the electric field $E_\varepsilon$ decreases as well. In particular, by taking the extreme case with
$d=\varepsilon/\eta$ and $0<\eta<1$, one recovers the configuration of the periodic structure in the homogenization regime (H1).
A straightforward asymptotic expansion of $\eqref{eq-p_H1}$ and $\eqref{eq-q_H1}$ for $p$ and $q$ leads to the following Lemma.
\begin{lem}\label{lem-pq_small_H2}
If $\varepsilon \ll 1$, $k=\varepsilon^\sigma$, and  $d=\varepsilon/\eta$, then 
\begin{equation*}
\dfrac{1}{p} \cdot \dfrac{1}{k\sin k} = \dfrac{1}{ 2 \alpha}  (1+O(\varepsilon^{\sigma})),
\end{equation*}
and
\begin{equation*}
\dfrac{1}{q} =\frac{i\cos\theta}{\eta\alpha}\varepsilon^\sigma \left(1+O(\varepsilon^{\sigma}) \right).
\end{equation*}
\end{lem}
Then a similar calculation as in Theorem \ref{thm-u_slit} and  \ref{thm-E_slit} yields that
$$ u = 1+ O(\varepsilon^{\sigma}) \quad\mbox{and}\quad  E_{\varepsilon,1} = O(1) $$
in the slits. That is, no enhancement is gained for such configuration.

\section{Conclusion}{\label{sec-conclusion}}
In this series of two papers, we have investigated the field enhancement and anomalous diffraction 
for electromagnetic wave scattering by a periodic array of perfect conducting subwavelength slits. The quantitative analysis
of the wave field is presented in both the diffraction regime and the homogenization regime. It is demonstrated
that the field enhancement in the diffraction regime is mainly attributed to scattering resonances.
Such enhancement becomes weaker if the resonant frequency is close to the Rayleigh anomaly.
In the homogenization regimes, the field enhancement can be attributed to certain non-resonant phenomenon.
In addition, surface plasmonic effect mimicking that of plasmonic metal exists for the periodic structure
with small period, and almost total transmission can be obtained for certain incident plane waves.

Based on the studies for the single slit case in \cite{lin_zhang16} and the periodic case in this series,
the mechanism for the field enhancement and anomalous diffraction for  perfect conducting subwavelength slits is now clearly understood.
Along this line of research, we will explore the field enhancement and anomalous diffraction (or transmission) for a single narrow slit and an array of slits with plasmonic metals.
Other than the mechanisms that are already known to occur for perfect conductors, it is expected that additional
enhancement mechanisms, including surface plasmonic resonances, will be present.
This will be reported in forthcoming papers.

\appendix
\section{Proof of Lemma \ref{lem-optK} for (H1)}
We prove Lemma \ref{lem-optK} for the homogenization regime (H1) in this section.
Let $\Omega_1= (0, \frac{1}{\eta})\times (0, \infty)$, $\Omega_2= (0, 1)\times (0, -\infty)$, 
and $\Omega_{1,N}= (0, \frac{1}{\eta})\times (0, N)$, 
$\Omega_{2,N}= (0, 1)\times (0, -N)$. We first introduce two Green's functions for the domain $\Omega_1$ and $\Omega_2$ respectively. 

For $x, y \in \Omega_1$, we define 
\begin{align*}
G_1(x, y)&= -\sum_{n=1}^{\infty} \frac{1}{2n\pi \eta} \left(e^{-2n\pi \eta|x_2-y_2|}+ e^{-2n\pi \eta|x_2+y_2|}\right) \left(\cos{2n\pi \eta x_1}\cos{2n\pi \eta y_1} + \sin{2n\pi \eta x_1}\sin{2n\pi \eta y_1}\right)\\
&= -\sum_{n=1}^{\infty} \frac{1}{2n\pi \eta} \left(e^{-2n\pi \eta|x_2-y_2|}+ e^{-2n\pi \eta|x_2+y_2|}\right)\cos{2n\pi \eta(x_1-y_1)}.
\end{align*} 
It is clear that $G_1$ satisfies the following equations:
\begin{equation*} 
\left\{
\begin{array}{lllll}
\vspace*{0.2cm}
& \Delta_x G_1(x,y) = \delta (x-y), \quad \mbox{ for }\, x \in \Omega_1,\\
\vspace*{0.2cm}
& \dfrac{\partial G_1(x,y)}{\partial x_2}= 0, \quad \mbox{for }\,\,x_2=0,\\
\vspace*{0.2cm}
& G_1(0,x_2,y) = G_1(1/\eta,x_2,y), \\
\vspace*{0.2cm}
& \displaystyle\int_{0}^{\frac{1}{\eta}} G_1(x_1, 0,y) dx_1  =0,\\
\vspace*{0.2cm}
&G_1(\cdot, y) \to 0 \,\,\mbox{as $x_2 \to \infty$ and satisfies the outgoing radiation condition (\ref{eq-rad_cond})}.
\end{array}
\right.
\end{equation*}
Moreover, when both $x, y$ are restricted to the boundary $\{(x_1, x_2): x_1 \in (0, \frac{1}{\eta}), x_2 =0\}$, we have
\[
G_1(x_1, 0, y_1, 0) = -\sum_{n=1}^{\infty} \frac{1}{n\pi} \cos{2n\pi \eta(x_1-y_1)}=\frac{1}{\pi}\ln|2\sin{ \pi \eta (x_1-y_1)}|.
\]

For $x, y \in \Omega_2$,
we define 
\[
G_2(x, y)= -\sum_{n=1}^{\infty} \frac{1}{2n\pi} \left(e^{-2n\pi |x_2-y_2|}+ e^{-2n\pi|x_2+y_2|}\right)\cos{n\pi x_1}\cos{n\pi y_1}.
\]
Then $G_2$ solves the following equations:
\begin{equation*} 
\left\{
\begin{array}{llllll}
\vspace*{0.2cm}
&\Delta_x G_2(x,y) = \delta (x-y), \\
\vspace*{0.2cm}
&\dfrac{\partial G_2(x,y)}{\partial x_2}= 0, \quad \mbox{for }\,\,x_2=0,\\
\vspace*{0.2cm}
&\dfrac{\partial G_2(x,y)}{\partial x_1}= 0, \quad \mbox{for }\,\,x_1=0 \; \mbox{and} \;  x_1=1,\\
\vspace*{0.2cm}
&\displaystyle\int_{0}^{1} G_2(x_1, 0,y) dx_1  =0,\\
\vspace*{0.2cm}
& G_2(\cdot, y) \to 0 \,\,\mbox{as $x_2 \to -\infty$ and satisfies the outgoing radiation condition (\ref{eq-rad_cond})}.
\end{array}
\right.
\end{equation*}
Moreover,  when both $x, y$ are restricted to the boundary $\{(x_1, x_2): x_1 \in (0, 1), x_2 =0\}$, we have
\begin{eqnarray*}
G_2(x_1, 0, y_1, 0) &=& -\sum_{n=1}^{\infty} \frac{1}{n\pi} \left(\cos{n\pi(x_1-y_1)} +\cos{n\pi(x_1+y_1)}\right)  \\
 &=& \frac{1}{\pi}\ln\abs{4\sin{ \frac{\pi(x_1-y_1)}{2}}\sin{  \frac{\pi(x_1+y_1)}{2}}}.
\end{eqnarray*}

Recall that $V_1$ is the space of distributions in $H^{-\frac{1}{2}}(\mathbf{R})$ whose support is contained in $[0, 1]$, or distributions  defined in the interval $(0,1)$ whose zero extension to the whole line belongs to $H^{-\frac{1}{2}}(\mathbf{R})$. For any $\psi \in V_1$, 
we define two functions
\begin{eqnarray}\label{eq-u1u2}
u_1&=& K_1 \psi (x_1, x_2)= \int_{0}^{\frac{1}{\eta}} G_1(x_1, x_2, y_1, 0) \psi(y_1) dy_1=\int_{0}^1 G_1(x_1, x_2, y_1, 0) \psi(y_1) dy_1,\label{eq-u1} \\
u_2&=&K_2 \psi (x_1, x_2)= \int_{0}^1 G_2(x_1, x_2, y_1, 0) \psi(y_1) dy_1. \label{eq-u2}
\end{eqnarray}
By the Green's identity, one can show that $u_1$ and $u_2$ is the unique solution to the problem
\begin{equation*} 
\left\{
\begin{array}{lllll}
\vspace*{0.2cm}
& \Delta u_1(x) =0, \quad \mbox{ for }\, x \in \Omega_1,\\
\vspace*{0.2cm}
& \dfrac{\partial u_1(x)}{\partial x_2}= \psi, \quad \mbox{for }\,\,x_2=0,\\
\vspace*{0.2cm}
& u_1(0,x_2)=u_1(1/\eta,x_2),  \\
& \displaystyle\int_{0}^{\frac{1}{\eta}} u_1(x_1, 0) dx_1  =0.\\
\vspace*{0.2cm}
& u_1 \to 0 \,\,\mbox{as $x_2 \to \infty$},
\end{array}
\right.
\quad \mbox{and} \quad
\left\{
\begin{array}{lllll}
\vspace*{0.2cm}
& \Delta u_2(x) =0, \quad \mbox{ for }\, x \in \Omega_2,\\
\vspace*{0.2cm}
& -\dfrac{\partial u_2(x)}{\partial x_2}= \psi, \quad \mbox{for }\,\,x_2=0,\\
\vspace*{0.2cm}
& u_2(0,x_2)=u_2(1,x_2),  \\
& \displaystyle\int_{0}^{1} u_2(x_1, 0) dx_1  =0.\\
\vspace*{0.2cm}
& u_2 \to 0 \,\,\mbox{as $x_2 \to -\infty$},
\end{array}
\right.
\end{equation*}
respectively.

Let us define the following two operators associated with the trace of the functions $u_1, u_2$:
\begin{align*}
K_{1,0} \psi (x_1)&= \int_{0}^1 G_1(x_1, 0, y_1, 0) \psi(y_1) dy_1, \\
K_{2,0} \psi (x_1)&= \int_{0}^1 G_2(x_1, 0, y_1, 0) \psi(y_1) dy_1.
\end{align*}
Let $K_0= K_{1, 0} + K_{2, 0}$. By analyzing the singularities in the kernel of the two operators $K_{1, 0}, K_{2, 0}$, and using the argument in \cite{eric10}, 
it follows that $K_0: V_1 \to (V_1)^*=V_2$ is bounded. Moreover, $K_0^*=K_0$, where 
$K_0^*$ is the dual operator (see \cite{conway}) of $K_0$. 

We show that 
\begin{lem} \label{lem-91}
$K_0$ is invertible from $V_1$ to $(V_1)^*$ and its inverse is bounded. 
\end{lem}

To establish the above result, we first prove the following identity.

\begin{lem}
For any $\psi \in V_1$, we have
\begin{align*}
\langle K_0 \psi, \psi \rangle&=-\int_{\Omega_{1}} |\nabla u_1|^2 dx_1dx_2- \int_{\Omega_{2}} |\nabla u_2|^2 dx_1dx_2,
\end{align*}
where $u_1$ and $u_2$ are defined in \eqref{eq-u1} and  \eqref{eq-u2} .
\end{lem}
\noindent\textbf{Proof}.
Note that both $u_1$ and $u_2$ are harmonic functions and can be expanded as
\begin{align*}
u_1 &= \sum_{n> 0} \left(  a_{n,1} \sin{2n\pi \eta x_1} + b_{n,1} \cos{2n\pi \eta x_1}\right) e^{-2n\pi \eta x_2},\\
u_2 &= \sum_{n> 0}  b_{n,2} \cos{n\pi x_1} e^{2n\pi x_2}.
\end{align*}
for some constants $a_{n,1}, b_{n,1}, b_{n,2}$. 

On the other hand, from the boundary conditions
\[
\frac{\partial u_1 (x_1, 0)}{ \partial x_2} = -\frac{\partial u_2 (x_1, 0)}{ \partial x_2}=\psi,
\]
it follows that
\begin{align*}
\int_{\Omega_{1}} |\nabla u_1|^2 dx_1dx_2 &=\lim_{N\to \infty}\int_{\Omega_{1,N}} |\nabla u_1|^2 dx_1dx_2=  
\lim_{N\to \infty}\int_{\partial \Omega_{1,N}} u_1 \frac{\partial u_1}{\partial \nu} d\sigma \\
&= -\int_{0}^1 u_1 (x_1, 0) \frac{\partial u_1(x_1, 0)}{\partial x_2} dx_1 + \lim_{N\to \infty} \int_{0}^{\frac{1}{\eta}}u_1(x_1, N) \frac{\partial u_1(x_1, N)}{\partial x_2} dx_1 \\
&= -\langle K_{1, 0} \psi, \psi \rangle - \lim_{N\to \infty}\sum_{n=1}^{\infty} n \pi \eta e^{-4n\pi \eta N} \left( |a_{n,1}|^2 + |b_{n,1}|^2 \right) \\
&= -\langle K_{1, 0} \psi, \psi \rangle. 
\end{align*}

Similarly for $u_2$, we have 
\[
\int_{\Omega_{2}} |\nabla u_2|^2 dx_1dx_2 =-\langle K_{2, 0} \psi, \psi \rangle.
\]

The lemma follows.\\

Based on the above identity, we can show that 
\begin{lem} \label{lem-99}
There exists $C>0$ such that for all $\psi \in V_1$
\[
\|  K_0 \psi \|_{V_2}  \geq C \|\psi\|_{V_1}. 
\]
\end{lem}
\noindent\textbf{Proof}.
We consider $u_2$ restricted to the domain $\Omega_{2, 1}$. We have 
$$\int_{\Omega_{2,1}} u_2 dx_1dx_2 =0.$$
By Poincare's inequality, there exists a constant $C_1$ such that 
$$
\|u_2\|_{H^1(\Omega_{2, 1})} \leq C_1 \| \nabla u_2\|_{L^2(\Omega_{2, 1})} \leq C_1 \sqrt{\|  K \psi \|_{V_2} \cdot \| \psi\|_{V_1}}.
$$ 
On the other hand, note that $\psi = -\dfrac{\partial u_2 (x_1, 0)}{ \partial x_2}$. By the trace theorem, we have
$$
\|\psi\|_{V_1} \leq C_2  \|u_2\|_{H^1(\Omega_{2, 1})}
$$
for some constant $C_2$.
It follows that 
$$
\|\psi\|_{V_1} \leq C_1^2 C_2^2 \|  K \psi \|_{V_2}.
$$
This proves the lemma. \\

\noindent\textbf{Proof of Lemma \ref{lem-91}}. From Lemma \ref{lem-99}, we can conclude that the map $K_0: V_1 \to V_1^*$ is injective. This also shows that $K_0^*$ is also injective (since $K_0^*= K_0$). As a result, $K_0 (V_1)$ is dense in $(V_1)^*$. But Lemma \ref{lem-99} also implies that $K_0 (V_1)$ is closed in $V_2$. Therefore, 
$K_0(V_1)= (V_1)^*$ and consequently $K_0$ is has a bounded inverse $K_0^{-1}$ by the open mapping theorem.\\

\noindent\textbf{Proof of Lemma Lemma \ref{lem-optK} for the case H1}. For any $\psi \in V_1$, note that
$$
K \psi (X) = K_0 \psi - \frac{3\ln 2}{\pi} \langle \psi, 1 \rangle 1  + \dfrac{\kappa \eta }{ \sqrt{k^2-\kappa^2}} \int_0^1(X-Y)\psi(Y)dY.
$$
A direct calculation yields 
$$
\langle K \psi , \psi \rangle = \langle K_0 \psi , \psi \rangle - \frac{3\ln 2}{\pi} \abs{\langle \psi, 1\rangle}^2 < \langle K_0 \psi , \psi \rangle
$$
Therefore, using Lemma \ref{lem-99}, we can show that 
$$
\|K \psi \|_{V_1^*} \geq C\|\psi\|_{V_1}
$$
for some constant $C$. Similar to the proof of Lemma \ref{lem-91}, we can conclude that $K$ is invertible from $V_1$ to $V_1^*$ and its inverse is also bounded. 

To calculate $\alpha(k,\kappa):=\langle K^{-1} 1, 1 \rangle$. Let
$\psi_0 = K^{-1} 1$. Then $\psi_0$ depends on $k$ and $\kappa$ and we have
$$
\alpha(k,\kappa) = \langle \psi_0, K \psi_0\rangle =  \langle K_0 \psi_0 , \psi_0 \rangle- \frac{3\ln 2}{\pi} \abs{\langle \psi, 1\rangle}^2 <0.
$$
It is obvious that $\alpha(k,\kappa)$ is a real number. This completes the proof of Lemma \ref{lem-optK}.

\bibliography{references}

\end{document}